\documentclass[12pt]{article}
\usepackage[a4paper,margin=31mm]{geometry}
\usepackage{amsmath,amssymb,amsthm,mathtools}
\usepackage{microtype}
\usepackage[hidelinks]{hyperref}
\usepackage{enumitem}
\usepackage{xcolor}
\hypersetup{
 pdftitle={Growth of Sobolev Norms for Linear Schrodinger Equations with Denjoy--Carleman Potentials},
 pdfauthor={Yuan Wu, Huining Xue},
 pdfsubject={Quantitative finite normal forms and Sobolev norm growth},
 pdfkeywords={linear Schrodinger equation, Sobolev norms, Denjoy--Carleman class, normal form}
}

\newtheorem{theorem}{Theorem}[section]
\newtheorem{proposition}[theorem]{Proposition}
\newtheorem{lemma}[theorem]{Lemma}
\newtheorem{corollary}[theorem]{Corollary}
\theoremstyle{definition}
\newtheorem{definition}[theorem]{Definition}
\theoremstyle{remark}
\newtheorem{remark}[theorem]{Remark}

\newcommand{\T}{\mathbb T}
\newcommand{\Z}{\mathbb Z}

\newcommand{\R}{\mathbb R}

\newcommand{\la}{\langle}
\newcommand{\ra}{\rangle}

\title{Growth of Sobolev Norms for Linear Schr\"odinger\\
Equations with Denjoy--Carleman Potentials}
\author{Yuan Wu \and Huining Xue}

\begin{document}
\maketitle

\begin{abstract}
We consider the linear Schr\"odinger equation on a torus of arbitrary dimension with a
real time-dependent potential satisfying a mixed space--time Denjoy--Carleman condition.
More precisely, for an increasing concave weight
\(\omega(r)\asymp \log^\sigma(e+r)\), \(\sigma> 1\),we assume
\[
 \sup_{t\in\mathbb R}|\partial_t^q\widehat V(k,t)|
 \le Ae^{Bq^{p}}e^{-\rho_{0}\omega(|k|)},\qquad
 p=\frac{\sigma}{\sigma-1}.
\]
The corresponding spatial class is strictly larger than every Gevrey class, while the
moment factor gives the quantitative ultradifferentiable control in time needed by the
finite normal form.  We adapt the method of Delort and Fang--Zhang to this non-power
Fourier weight.  The essential new point is a quantitative operator calculus whose
majorant is built from the moment sequence and an auxiliary factorial that closes the
exact Leibniz convolutions.
For every $s>0$ we prove
\[
 \|u(t)\|_{H^s}\le C_s
 \exp\{C_s[\log(2+|t|)]^{\frac{1}{\sigma}}\}
 \|u(0)\|_{H^s}.
\]
The estimate is valid for every \(t\in\mathbb R\), holds in every space dimension, requires neither
smallness, quasi-periodicity, nor a transversality condition.
\end{abstract}



\section{Introduction}

Let $d\ge1$ and consider
\begin{equation}\label{eq:sch}
 i\partial_tu-\Delta u+V(x,t)u=0,
 \qquad (x,t)\in\T^d\times\R,
 \qquad u(0)=u_0.
\end{equation}
We use the convention
\(\T^d=(\mathbb R/2\pi\mathbb Z)^d\).  Thus
\[
 e_n(x)=(2\pi)^{-d/2}e^{in\cdot x},\qquad
 \Delta e_n=-|n|^2e_n,\qquad n\in\mathbb Z^d.
\]
We denote by \(\Pi_nu=\langle u,e_n\rangle e_n\) the orthogonal projector onto the
\(n\)-th Fourier mode and normalize
\(\widehat f(n)=(2\pi)^{-d}\int_{\mathbb T^d}f(x)e^{-in\cdot x}\,dx\).
Throughout, \(\mathbb N_0=\{0,1,2,\ldots\}\).  The letter \(C\) denotes a positive
constant which may change from line to line; decorated constants retain the dependencies
specified in their statements.  None of these constants depends on the terminal
normal-form order unless this is explicitly indicated.
The potential is assumed real.  For a merely smooth bounded potential, the standard
energy estimate gives polynomial growth of high Sobolev norms; the work of Bourgain
\cite{BourgainSmooth} and Delort \cite{Delort} gives bounds $O(\la t\ra^\varepsilon)$
for every $\varepsilon>0$.  Wang \cite{Wang} proved a logarithmic bound in one dimension
for analytic or Gevrey potentials, and Fang and Zhang \cite{FangZhang} extended the
Gevrey estimate to tori of arbitrary dimension.  Their proof is a quantitative finite
normal form: one block-diagonalizes the equation to an arbitrarily high smoothing order,
keeps precise track of the factorial constants, and then optimizes the order.

Earlier adiabatic approaches for abstract Hamiltonians with increasing spectral gaps were
developed by Nenciu \cite{Nenciu} and Barbaroux--Joye \cite{BarbarouxJoye}.
Related abstract normal-form approaches were developed by Maspero--Robert
\cite{MasperoRobert} and Bambusi--Gr\'ebert--Maspero--Robert \cite{BGMR}.  The former
also obtains logarithmic bounds under analytic time dependence in an abstract spectral-gap
setting.  The point here is different: on an arbitrary-dimensional torus the Laplace
spectrum has large multiplicities and no increasing consecutive gaps, so the argument must
combine the Bourgain--Delort lattice partition with the precise Fourier weight in both
space and time.

More recent work places the problem in a broader dynamical context.  Bambusi--Langella--
Montalto \cite{BLM} treat unbounded time-dependent perturbations of the Laplacian on flat
tori and obtain subpolynomial bounds, whereas Thomann \cite{Thomann} constructs an explicit
linear time-dependent Schr\"odinger model with optimal Sobolev growth.  Maspero
\cite{MasperoDispersive} gives a complementary dispersive mechanism producing unbounded
Sobolev orbits for periodic perturbations in constant-gap settings.  These results address
different regularity and dynamical regimes; the present contribution is the quantitative
upper bound obtained by passing beyond the Gevrey scale for a multiplication potential.

The purpose of this paper is to carry out this program for the larger Denjoy--Carleman
class singled out in \cite{ChenLiYuan}; for the Fourier-analytic theory of general
ultradifferentiable classes we refer to Braun--Meise--Taylor
\cite{BraunMeiseTaylor}.  Its characteristic Fourier decay is
\begin{equation}\label{eq:rawweight}
 |\widehat f(k)|\le C\exp\!\left[-c\log^\sigma(e+|k|)\right],
 \qquad \sigma>1.
\end{equation}
It is slower than $e^{-c|k|^\epsilon}$ for every $\epsilon>0$, and hence the class is
strictly larger than every Gevrey class.  It is nevertheless faster than every negative
power, which is exactly what permits a finite normal form of arbitrary order.

The contribution of the paper is deliberately quantitative rather than structural.  We do
not introduce a new spectral clustering theorem: we use the Bourgain--Delort partition of
the spectrum of the Laplacian.  Nor do we change the normal-form mechanism of
Delort--Fang--Zhang.  What is new is that every estimate in that mechanism is rebuilt for
the non-power weight in \eqref{eq:rawweight}.  In particular, we prove the weighted
convolution estimates, the moment inequalities needed under repeated composition, and the
uniform bound for all terms of a finite normal form.  These estimates show precisely that
the constant at smoothing order $N$ is
\[
 \exp\bigl(CN^{p}\bigr).
\]
The optimization of this constant produces \eqref{eq:mainbound}.  Thus the result extends
the admissible potentials beyond the Gevrey scale, while the dynamical argument remains the
finite-order block diagonalization developed in the cited works.

Let us also explain why no small-divisor hypothesis occurs.  A Floquet reduction of a
quasi-periodic potential $V(x,\nu t)$ would lead to denominators
$|n|^2-|n'|^2+\ell\cdot\nu$, which can be arbitrarily small and therefore normally require
Diophantine estimates and parameter exclusion; Bourgain's quasi-periodic analysis
\cite{BourgainQP} is a fundamental example of this space--time approach.  We do not Fourier transform in time and
do not solve that Floquet homological equation.  Time is retained as an external variable.
At every normal-form step we solve only
\[
 [Q,\Delta]=-A,
 \qquad
 Q_{nn'}=-\frac{A_{nn'}}{|n|^2-|n'|^2},
\]
and only on the matrix entries for which the Bourgain--Delort decomposition gives
\[
 \bigl||n|^2-|n'|^2\bigr|
 \ge c\bigl(1+\max(|n|,|n'|)\bigr)^\delta.
\]
The complementary off-block entries are not divided by a small number.  Cluster separation
then forces $n$ and $n'$ to be far apart, and the spatial Denjoy--Carleman Fourier decay
makes those entries arbitrarily smoothing.  Terms produced by differentiating $Q(t)$ are
also not inverted; they enter the next, lower-order stage of a \emph{finite} normal form.
This dichotomy---divide only across a known spectral gap, and treat spatially separated
entries as smoothing---is the precise reason that neither a Diophantine condition nor a
transversality assumption is needed.

One cannot obtain the correct result by formally replacing a Gevrey exponent in
\cite{FangZhang}.  The constants in an $N$-th order smoothing estimate are now governed by
the moments of the weight in \eqref{eq:rawweight}.  Those moments grow as
\[
 \exp\big(CN^{p}\big),
\]
instead of a fixed power of $N!$.  Optimizing the resulting high-order estimate yields an
exponential of a power of $\log\log |t|$, rather than a power of $\log |t|$.

\subsection{Precise hypotheses and main result}

For $q\in\mathbb N_0=\mathbb{N}\setminus \{0\}$ and $ \sigma > 1, p=\frac{\sigma}{\sigma-1}$, we assume \(V\in C^\infty(\T^d\times\mathbb R)\) and impose the following mixed
Denjoy--Carleman condition: there are
$A,B,\rho_{0}>0$ such that
\begin{equation}\label{eq:DC}
 \sup_{t\in\R}|\partial_t^q\widehat V(k,t)|
 \le A e^{Bq^{p}}e^{-\rho_{0}\omega(|k|)},
 \qquad q\in\mathbb N_0,\quad k\in\Z^d.
\end{equation}
Here $\widehat V(k,t)$ denotes the Fourier coefficient in $x$.

\begin{theorem}[Denjoy--Carleman Sobolev bound]\label{thm:main}
Assume that $V$ is a real-valued function and satisfies \eqref{eq:DC}. Then for any $s>0$ and $u_0\in H^s(\T^d)$, the equation \eqref{eq:sch} has a unique global solution $ u(t) \in C(\R;H^s(\T^d))$.  Particularly, there are constants
$C_s,c_s>0$, depending on $s,d,\sigma $, and the constants $A,B,\rho_{0}$, such that
\begin{equation}\label{eq:mainbound}
 \|u(t)\|_{H^s}
 \le C_s\exp\{c_s\log(1+|t|)^{\frac{1}{\sigma}}\}
 \|u_0\|_{H^s},\qquad t\in\R.
\end{equation}
\end{theorem}

\begin{remark}
The transversality assumption in \cite{ChenLiYuan} is needed there for a large-deviation
argument and Anderson localization.  It plays no role in the present energy-growth problem.
Likewise, quasi-periodicity and arithmetic properties of a frequency are unnecessary.
\end{remark}

\section{The weight and its moment sequence}
 For $ 0<\rho <\rho_0$, $M\ge1$ and $\kappa>0$, set
\begin{eqnarray}\nonumber
 && \phi(r)=[\log(R_\sigma+r)]^\sigma-(\log R_\sigma)^\sigma, \quad R_\sigma=\max(e,e^{\sigma-1}),\\
\nonumber&& W_b=M^b(b!)^2e^{\kappa b^p},\quad b\in\mathbb{N}.
\end{eqnarray}
Then we have
\begin{lemma}[Properties of the weights]
For $r,s\ge0$,
\begin{equation}\label{eq:subadd}
 \phi(r+s)\le\phi(r)+\phi(s),\qquad
 |\phi(r)-[\log(e+r)]^\sigma|\le C_\sigma.
\end{equation}
and, for $b\ge0$ and $\eta>0$,
\begin{equation}\label{eq:moment}
 \sup_{r\ge1}r^b e^{-\eta(\log r)^\sigma}
 =\exp\bigl(D_\sigma\eta^{-1/(\sigma-1)}b^p\bigr),\quad
 D_\sigma=(\sigma-1)\sigma^{-\sigma/(\sigma-1)}.
\end{equation}
Furthermore, we have
$$\sup_{r\ge0}r^b e^{-\eta[\log(e+r)]^\sigma}\leq \exp\bigl(D_\sigma\eta^{-1/(\sigma-1)}b^p\bigr). $$
\end{lemma}
\begin{proof}
Direct differentiation gives
\[
\phi'(r)
=
\frac{\sigma[\log(R_\sigma+r)]^{\sigma-1}}
{R_\sigma+r}>0
\]
and
\[
\phi''(r)
=
\frac{\sigma[\log(R_\sigma+r)]^{\sigma-2}}
{(R_\sigma+r)^2}
\left(\sigma-1-\log(R_\sigma+r)\right)\leq 0.
\]
Hence $\phi$ is increasing and concave, and we also have $\phi(0)=0$.

If $r+s>0$, concavity implies
\[
\begin{aligned}
\phi(r)
&=
\phi\left(
\frac{r}{r+s}(r+s)+\frac{s}{r+s}\,0
\right)\\
&\ge\frac{r}{r+s}\phi(r+s).
\end{aligned}
\]
Similarly,
\[
\phi(s)\ge\frac{s}{r+s}\phi(r+s).
\]
Then one has
\[
\phi(r+s)\le\phi(r)+\phi(s).
\]
The case $r=s=0$ is immediate.

Next, since $R_\sigma\ge e$,
\[
[\log(R_\sigma+r)]^\sigma-[\log(e+r)]^\sigma
=
\int_e^{R_\sigma}
\frac{\sigma[\log(r+t)]^{\sigma-1}}{r+t}\,dt.
\]
The function
\[
x\longmapsto\frac{(\log x)^{\sigma-1}}{x}
\]
is continuous on $[e,\infty)$ and tends to zero as $x\to\infty$.
Consequently,
\[
L_\sigma:=
\sup_{x\ge e}\frac{(\log x)^{\sigma-1}}{x}<\infty,
\]
and
\[
0\le
[\log(R_\sigma+r)]^\sigma-[\log(e+r)]^\sigma
\le\sigma(R_\sigma-e)L_\sigma.
\]
Therefore
\[
\left|\phi(r)-[\log(e+r)]^\sigma\right|
\le
(\log R_\sigma)^\sigma
+\sigma(R_\sigma-e)L_\sigma
=:C_\sigma.
\]

To prove the moment identity, set $y=\log r$ for $r\ge1$.
Then
\[
r^b e^{-\eta(\log r)^\sigma}
=
e^{by-\eta y^\sigma},
\qquad y\ge0.
\]
For $b=0$, the exponent attains its maximum $0$ at $y=0$,
so the asserted identity holds.

Suppose now that $b>0$. Since
\[
\frac{d}{dy}(by-\eta y^\sigma)
=
b-\eta\sigma y^{\sigma-1},
\]
the unique maximizer is
\[
y_*=
\left(\frac{b}{\eta\sigma}\right)^{1/(\sigma-1)}.
\]
Indeed, the derivative is positive for $y<y_*$ and negative for
$y>y_*$. Then we obtain
\[
\sup_{r\ge1}r^b e^{-\eta(\log r)^\sigma}
=\exp\left(by_*-\eta y_*^\sigma\right)=
\exp\!\left(
D_\sigma\eta^{-1/(\sigma-1)}b^p
\right).
\]

Finally, for $r\ge0$ and $b\ge0$,
\[
r^b\le(e+r)^b.
\]
Setting $y=\log(e+r)\ge1$, we find
\[
r^b e^{-\eta[\log(e+r)]^\sigma}
\le e^{by-\eta y^\sigma}.
\]
It follows that
\[
\begin{aligned}
\sup_{r\ge0}r^b e^{-\eta[\log(e+r)]^\sigma}
\le\sup_{y\ge1}e^{by-\eta y^\sigma}\le\sup_{y\ge0}e^{by-\eta y^\sigma}=\exp\!\left(
D_\sigma\eta^{-1/(\sigma-1)}b^p
\right),
\end{aligned}
\]
which completes the proof.

\end{proof}

\section{Spectral blocks and a fixed majorant sequence}\label{sec:blocks}

\begin{proposition}[Bourgain--Delort partition]
\label{prop:clusters}
For $ 0<\sigma_{0}< 1/10$, there are \(0<\tau_{0}<\sigma_0\), a constant \(\gamma>0\), and a partition
\[
 \Z^d=\bigcup_{\alpha\in\mathfrak A}A_\alpha
\]
and representatives \(n_\alpha\in A_\alpha\) such that
\begin{align}
 n,n'\in A_\alpha&\Longrightarrow
 |n-n'|+\big||n|^2-|n'|^2\big|
 \le \gamma+\max(|n|,|n'|)^{\sigma_0},\label{eq:newinside}\\
 n\in A_\alpha,\ n'\in A_\beta,\ \alpha\ne\beta&\Longrightarrow
 |n-n'|+\big||n|^2-|n'|^2\big|
 \ge \max(|n|,|n'|)^{\tau_{0}} .
\label{eq:newsep}
\end{align}
Moreover \(\langle n\rangle\asymp\langle n_\alpha\rangle\) on \(A_\alpha\).
\end{proposition}

\begin{proof}
This is \cite[Lemma 3.2]{Delort}; the identical formulation used for the torus appears as
\cite[Proposition 2.1]{FangZhang}.  We only record the consequence needed for the
equivalent Sobolev norm.  Put \(R=\max(|n|,|n_\alpha|)\).  By
\eqref{eq:newinside}, \(|n-n_\alpha|\le\gamma+R^{\sigma_0}\).  Since
\(\sigma_0<1\), the right-hand side is at most \(R/2\) outside a fixed ball; hence
\(\min(\langle n\rangle,\langle n_\alpha\rangle)\ge cR\).  Enlarging the comparison
constant on the fixed ball proves the final assertion.
\end{proof}

Write \(\Pi_\alpha=\sum_{n\in A_\alpha}\Pi_n\) and set
\[
 \widetilde\Delta=-\sum_\alpha |n_\alpha|^2\Pi_\alpha,\qquad
 \|u\|_{\widetilde H^s}=\|(1-\widetilde\Delta)^{s/2}u\|_2.
\]
Then
\begin{equation}\label{eq:normequi}
 C_0^{-|s|}\|u\|_{H^s}\le \|u\|_{\widetilde H^s}
 \le C_0^{|s|}\|u\|_{H^s}.
\end{equation}
Indeed, both sides squared are diagonal sums, and the assertion follows term by term from
\(\langle n\rangle\asymp\langle n_\alpha\rangle\).

For \(n,n'\in\Z^d\), put
\[
 h(z)=\langle z\rangle^{-(d+2)},\quad r_{nn'}=\max(|n|,|n'|),\quad X_{nn'}=1+r_{n,n'}, \quad b_{+}=\max(b,0).
 \]
All operator classes below are defined with estimates for every
\(q,N\in\mathbb N_0\).  The normal form itself is finite, but no derivative cutoff is
built into the definitions.  This removes endpoint ambiguities when a time derivative or
a Laplacian commutator is taken.

\begin{definition}[Principal classes and the derivative flag]\label{def:L} For $j\in\mathbb{N}$ and $\delta\in\{0,1\}$, let $\mathcal{L}_\tau^{-j}(\delta)=\mathcal{L}_{\tau,\sigma}^{-j}(M,\kappa,\rho;\delta)$ be the class of smooth time-dependent operator families
$Q(t):C^\infty(\mathbb{T}^d)\to\mathcal D'(\mathbb{T}^d)$ for which there exists $B_Q\ge0$ such that, for all $q,n,n'$,
\begin{eqnarray}\label{eq:L}
 \sup_t\|\Pi_n\partial_t^qQ(t)\Pi_{n'}\|_{2\to2}
 &\le& B_Q W_{q+(j+\delta-1)_+}
 e^{-\rho\phi(|n-n'|)}h(n-n')X_{nn'}^{-j\tau}\\
 &&\times\mathbf1_{\{|n-n'|\le r_{nn'}/[10(1+j)]\}}.
\end{eqnarray}
The least such $B_Q$ is denoted by $ \|Q\|_{j,\delta}$. Once the parameters are fixed, $B_Q$ must be uniform in $q,n,n',t$.
\end{definition}

For $\widetilde{\mathcal{L}}_\tau^{-j}(\delta)$, replace the cutoff in~\eqref{eq:L} by
\[
 \widetilde\chi_j(n,n')
 =\mathbf1_{\{|n-n'|\le r_{nn'}/[10(2+j)]\}}.
\]
The class $\overline{\mathcal{L}}_\tau^{-j}(\delta)$ additionally requires the matrix support to lie in
\[
 \bigl||n|^2-|n'|^2\bigr|>\frac14(|n|+|n'|)^{\tau_0}.
\]
The seminorm is again the least constant in the corresponding matrix estimate. An overline denotes the large-gap subclass, distinguishing it from the ordinary class $\mathcal{L}$.

\begin{definition}[Block-diagonal and off-block subclasses]\label{2}
For the narrower-cutoff or large-gap classes above, the subscript $\mathrm D$ means
$\widetilde\Pi_\alpha Q\widetilde\Pi_\beta=0$ ($\alpha\ne\beta$),
and the subscript $\mathrm{ND}$ means $\widetilde\Pi_\alpha Q\widetilde\Pi_\alpha=0$.

Let
\[
 Q_{\mathrm D}=\sum_\alpha\widetilde\Pi_\alpha Q\widetilde\Pi_\alpha,
 \qquad Q_{\mathrm{ND}}=\sum_{\alpha\ne\beta}\widetilde\Pi_\alpha Q\widetilde\Pi_\beta.
\]
Block projections do not increase the corresponding seminorms. Throughout, $Q^*(t)$ denotes the $L^2$ adjoint at fixed time. Since the weights and cutoffs are symmetric in $n,n'$, taking the adjoint also preserves the seminorm.
\end{definition}

\begin{proposition}[Commutation of block-diagonal operators]\label{2.2}
If $Q$ belongs to $\widetilde{\mathcal{L}}_{\tau,D}^{-j}(\delta)$ or $\bar{\mathcal{L}}_{\tau,D}^{-j}(\delta)$, then $[Q,\widetilde\Delta]=0$.
\end{proposition}
\begin{proof} The operator $\widetilde\Delta$ is scalar on each cluster; the claim follows block by block.\end{proof}

\begin{proposition}[The homological equation]\label{2.3}
Fix $0<\tau\le\tau_0$. If $j\ge1$ and
$S\in\bar{\mathcal{L}}_{\tau,\mathrm{ND}}^{-(j-1)}(\delta)$, $\delta\in\{0,1\}$, define
\begin{equation}\label{eq:homological}
 \Pi_nQ\Pi_{n'}=
 \frac{\Pi_nS\Pi_{n'}}{|n'|^2-|n|^2}
\end{equation}
on the support of $S$, and set all other entries to zero. Then
\[
 [Q,\Delta]=-S,\qquad Q\in\mathcal{L}_\tau^{-j}(0),\qquad
 \| Q\|_{j,0}\le C\|S\|_{j-1,\delta}.
\]
If $S^*=S$, then $Q^*=-Q$, where $ Q^*$ denotes the adjoint of $ Q$.
\end{proposition}
\begin{proof}
The reciprocal of the large-gap denominator is at most $CX_{nn'}^{-\tau_0}$, while
$(j-2+\delta)_+\le j-1$. The narrower input cutoff is exactly the ordinary cutoff at output level $j$. Skew-adjointness follows because the denominator changes sign when $n,n'$ are interchanged.
\end{proof}

\begin{definition}[Smoothing remainders]\label{3}
Let $\mathcal{R}_j^{-\infty}=\mathcal{R}_{j,\sigma}^{-\infty}(M,\kappa,\tau)$ consist of smooth operator families satisfying
\begin{equation}\label{eq:R}
 \sup_t\|\Pi_n\partial_t^qR(t)\Pi_{n'}\|_{2\to2}
 \le B_R W_{N+j+q}h(n-n')X_{nn'}^{-\tau N},
 \qquad N,q\in\mathbb{N},
\end{equation}
and denote the least admissible $B_R$ by $|R|_j$. The same $B_R$ must work for all $N,q$. This definition uses the combined weight $W_{N+j+q}$.
\end{definition}

\begin{proposition}[One time derivative]\label{2.4}
For $j\ge1$ and $Q\in\mathcal{L}_\tau^{-j}(0)$,
\[
 [i\partial_t,Q]=i\partial_tQ\in\mathcal{L}_\tau^{-j}(1),\qquad
 \|[i\partial_t,Q]\|_{j,1}\le \|Q\|_{j,0}.
\]
\end{proposition}
Indeed, $q+1+(j-1)_{+}=q+j$. The case $j=0$ is excluded; the same shift formula cannot be applied directly at that endpoint.

\begin{proposition}[Sobolev mapping properties of principal operators]\label{2.5}
If $Q\in\mathcal{L}_\tau^{-j}(\delta)$, then, for $s\in\mathbb{R}$ and $q\in\mathbb{N}$,
\begin{equation}\label{eq:Lmap}
 \|\partial_t^qQ\|_{H^s\to H^{s+j\tau}}
 \le C_d \left(\frac{10}{9}\right)^{|s|}\| Q\|_{j,\delta}W_{q+(j+\delta-1)_+}.
\end{equation}
\end{proposition}

\begin{proof}
We use the notation
\[
\langle n\rangle=(1+|n|^2)^{1/2},
\qquad
r_{nn'}=\max(|n|,|n'|),
\qquad
X_{nn'}=1+r_{nn'},
\]
and
\[
h(k)=\langle k\rangle^{-(d+2)}.
\]
Fix $q\in\mathbb N_0$ and set
\[
A_q(t)=\partial_t^qQ(t),
\qquad
\Lambda_q=\|Q\|_{j,\delta}W_{q+(j+\delta-1)_+}.
\]
By the definition of $\mathcal L_\tau^{-j}(\delta)$,
\begin{equation}
\begin{split}
\|\Pi_nA_q(t)\Pi_{n'}\|_{L^2\to L^2}
\le{}&
\Lambda_q
e^{-\rho\phi(|n-n'|)}
h(n-n')X_{nn'}^{-j\tau}\times
\mathbf 1_{\{
|n-n'|\le r_{nn'}/[10(1+j)]
\}},
\end{split}
\end{equation}
uniformly in $t$.

On the support of the matrix estimate,
\[
|n-n'|
\le\frac{r_{nn'}}{10(1+j)}
\le\frac{r_{nn'}}{10}.
\]
The reverse triangle inequality therefore gives
\[
\min(|n|,|n'|)
\ge r_{nn'}-|n-n'|
\ge\frac9{10}r_{nn'}.
\]
Consequently,
\[
\begin{split}
\min(\langle n\rangle,\langle n'\rangle)
\ge
\left(1+\frac{81}{100}r_{nn'}^2\right)^{1/2}\ge
\frac9{10}(1+r_{nn'}^2)^{1/2}=
\frac9{10}\max(\langle n\rangle,\langle n'\rangle).
\end{split}
\]
Thus
\[
\frac9{10}
\le
\frac{\langle n\rangle}{\langle n'\rangle}
\le
\frac{10}{9}.
\]
For every $s\in\mathbb R$, it implies
\begin{equation}\label{251}
\left(
\frac{\langle n\rangle}{\langle n'\rangle}
\right)^s
\le
\left(\frac{10}{9}\right)^{|s|}.
\end{equation}

Since
\[
\langle n\rangle\le1+|n|\le X_{nn'},
\]
and $j\tau\ge0$, we have
\[
\langle n\rangle^{j\tau}X_{nn'}^{-j\tau}\le1.
\]
Combining with (\ref{251}) yields, on the matrix support,
\begin{equation}\label{2511}
\begin{split}
\langle n\rangle^{s+j\tau}
X_{nn'}^{-j\tau}
\langle n'\rangle^{-s}
&=
\left(
\frac{\langle n\rangle}{\langle n'\rangle}
\right)^s
\left(
\frac{\langle n\rangle}{X_{nn'}}
\right)^{j\tau}\le
\left(\frac{10}{9}\right)^{|s|}.
\end{split}
\end{equation}
Let $u$ be a trigonometric polynomial and define
\[
a_{n'}=\langle n'\rangle^s\|\Pi_{n'}u\|_{L^2}.
\]
By the definition of the Sobolev norm,
\[
\|a\|_{\ell^2(\mathbb Z^d)}=\|u\|_{H^s}.
\]
Using
\[
\Pi_nA_q(t)u
=
\sum_{n'}\Pi_nA_q(t)\Pi_{n'}u,
\]
the matrix estimate and (\ref{2511}) give
\begin{align*}
\langle n\rangle^{s+j\tau}
\|\Pi_nA_q(t)u\|_{L^2}
&\le
\sum_{n'}
\langle n\rangle^{s+j\tau}
\|\Pi_nA_q(t)\Pi_{n'}\|_{L^2\to L^2}
\|\Pi_{n'}u\|_{L^2}\\
&\le
\Lambda_q
\left(\frac{10}{9}\right)^{|s|}
\sum_{n'}h(n-n')a_{n'}\\
&=
\Lambda_q
\left(\frac{10}{9}\right)^{|s|}
(h*a)_n.
\end{align*}
Since $d+2>d$, we have $h\in\ell^1(\mathbb Z^d)$.
The discrete Young inequality gives
\[
\|h*a\|_{\ell^2}
\le
\|h\|_{\ell^1}\|a\|_{\ell^2}.
\]
Taking the $\ell^2$ norm in the preceding pointwise estimate, we obtain
\begin{align*}
\|A_q(t)u\|_{H^{s+j\tau}}
&\le
\Lambda_q
\left(\frac{10}{9}\right)^{|s|}
\|h*a\|_{\ell^2}\le
C_d
\left(\frac{10}{9}\right)^{|s|}
\Lambda_q\|u\|_{H^s}.
\end{align*}

Recall that trigonometric polynomials are dense in $H^s(\mathbb T^d)$. For $u\in H^s$, choose trigonometric polynomials $u_N$ converging to $u$ in $H^s$. The estimate above implies
\[
\|A_q(t)(u_N-u_M)\|_{H^{s+j\tau}}
\le
C_d\left(\frac{10}{9}\right)^{|s|}
\Lambda_q\|u_N-u_M\|_{H^s}.
\]
Hence $A_q(t)u_N$ converges in $H^{s+j\tau}$.
Its limit defines the unique bounded extension of $A_q(t)$,
with the same norm bound. Since all estimates are uniform in
$t$, taking the supremum over $t\in\mathbb R$ proves the assertion.
\end{proof}

\begin{proposition}[Sobolev mapping properties of remainders]\label{2.6}
For $R\in\mathcal{R}_j^{-\infty}$, $s\ge0$, and $q,m\in\mathbb{N}$, set
$N_0=m+\lceil s/\tau\rceil$ and $N_1=m+\lceil(s+2)/\tau\rceil$. Then
\begin{align}
 \|\partial_t^qR\|_{H^{-s}\to H^{\tau m}}
 &\le C_d|R|_jW_{N_0+j+q},\label{eq:Rmap}\\
 \|\partial_t^q[\Delta,R]\|_{H^{-s}\to H^{\tau m}}
 &\le C_d|R|_jW_{N_1+j+q}.\label{eq:Rcomm}
\end{align}
\end{proposition}
\begin{proof}
The product of the two Sobolev weights is at most $X_{nn'}^{s+\tau m}$.
The commutator contributes at most an additional factor $2X_{nn'}^2$. Choose the corresponding $N$ and use the summability of $h$.
\end{proof}

Fix $0<\tau\le1$. Set
\[
 b_\rho=D_\sigma\rho^{-1/(\sigma-1)}\tau^p,
\]
and choose in advance
\begin{equation}\label{eq:threshold}
 \kappa\ge2b_\rho,\qquad
 \log M\ge\max\left(0,\sup_{j\in\mathbb{N}}
 \left\{\tau\log(10(j+2))-\frac{\kappa p}{2}j^{p-1}\right\}\right).
\end{equation}
The supremum is finite. These conditions are sufficient, with no claim of optimality.

\begin{proposition}[Composition of two principal operators with $\delta=0$]\label{2.7}
Let $Q_i\in\mathcal{L}_\tau^{-j_i}(0)$, $j_i\ge0$, and $j=j_1+j_2$. Define
\[
 \Pi_n\mathcal{M}(Q_1,Q_2)\Pi_{n'}
 =\Pi_nQ_1Q_2\Pi_{n'}\,\widetilde\chi_j(n,n'),
 \qquad \mathcal{E}(Q_1,Q_2)=Q_1Q_2-\mathcal{M}(Q_1,Q_2).
\]
Then
\begin{equation}\label{eq:prod}
 \mathcal{M}(Q_1,Q_2)\in\widetilde{\mathcal{L}}_\tau^{-j}(0),\qquad
 \mathcal{E}(Q_1,Q_2)\in\mathcal{R}_j^{-\infty},
\end{equation}
and
\[
 \|\mathcal{M}(Q_1,Q_2)\|_{j,0}+|\mathcal{E}(Q_1,Q_2)|_j
 \le C\|Q_1\|_{j_1,0}\|Q_2\|_{j_2,0}.
\]
The constant is independent of $j_1,j_2$.
\end{proposition}

\begin{proof}
Set
\[
j=j_1+j_2,\qquad
K_i=\|Q_i\|_{j_i,0},\qquad
b_i=(j_i-1)_+,\qquad b=b_1+b_2.
\]
Then
\[
b\le (j-1)_+\le j.
\]

Leibniz's formula gives
\begin{equation}\label{eq:p27-leibniz}
\begin{split}
\Pi_n\partial_t^q(Q_1Q_2)\Pi_{n'}
={}&
\sum_{\ell=0}^q\binom q\ell
\sum_{k\in\mathbb Z^d}
\bigl(\Pi_n\partial_t^\ell Q_1\Pi_k\bigr)
\bigl(\Pi_k\partial_t^{q-\ell}Q_2\Pi_{n'}\bigr).
\end{split}
\end{equation}
By the definitions of the operator classes, after dropping the
input support indicators,
\begin{align}
\|\Pi_n\partial_t^\ell Q_1\Pi_k\|_{2\to2}
&\le
K_1W_{\ell+b_1}
e^{-\rho\phi(|n-k|)}
h(n-k)X_{nk}^{-j_1\tau},
\label{eq:p27-input1}\\
\|\Pi_k\partial_t^{q-\ell}Q_2\Pi_{n'}\|_{2\to2}
&\le
K_2W_{q-\ell+b_2}
e^{-\rho\phi(|k-n'|)}
h(k-n')X_{kn'}^{-j_2\tau}.
\label{eq:p27-input2}
\end{align}
On the support of $\widetilde\chi_j(n,n')$, we have
\[
|n-n'|\le\frac{r_{nn'}}{10(j+2)}.
\]
Set
\[
\varepsilon_j=\frac1{10(j+2)},\qquad r=r_{nn'}.
\]
The reverse triangle inequality yields
\[
\min(|n|,|n'|)\ge(1-\varepsilon_j)r.
\]
Hence, for every intermediate frequency $k$,
\[
X_{nk}\ge1+|n|
\ge(1-\varepsilon_j)X_{nn'},
\]
and likewise
\[
X_{kn'}\ge(1-\varepsilon_j)X_{nn'}.
\]
It follows that
\[
X_{nk}^{-j_1\tau}X_{kn'}^{-j_2\tau}
\le
(1-\varepsilon_j)^{-j\tau}X_{nn'}^{-j\tau}.
\]
The prefactor is bounded independently of $j$, because
\[
-j\tau\log(1-\varepsilon_j)
\le
\frac{j\tau\varepsilon_j}{1-\varepsilon_j}
=
\frac{j\tau}{10j+19}
\le\frac{\tau}{10}.
\]
Thus
\begin{equation}\label{eq:p27-radial}
X_{nk}^{-j_1\tau}X_{kn'}^{-j_2\tau}
\le e^{\tau/10}X_{nn'}^{-j\tau}.
\end{equation}
Combining \eqref{eq:p27-leibniz}--\eqref{eq:p27-radial}
with the temporal convolution estimate (\ref{t1}) and spatial convolution estimates (\ref{s1}) gives
\begin{equation}\label{eq:p27-main}
\begin{split}
&\|\Pi_n\partial_t^q\mathcal M(Q_1,Q_2)\Pi_{n'}\|_{2\to2}\le
CK_1K_2W_{q+b}
e^{-\rho\phi(|n-n'|)}
h(n-n')X_{nn'}^{-j\tau}
\widetilde\chi_j(n,n').
\end{split}
\end{equation}
Since $b\le(j-1)_+$ and $W_m$ is nondecreasing,
\[
W_{q+b}\le W_{q+(j-1)_+}.
\]
Therefore
\[
\mathcal M(Q_1,Q_2)
\in\widetilde{\mathcal L}_\tau^{-j}(0),
\qquad
\|\mathcal M(Q_1,Q_2)\|_{j,0}\le CK_1K_2.
\]

For the remainder, use
\[
X_{nk}^{-j_1\tau}X_{kn'}^{-j_2\tau}\le1.
\]
The convolution estimates then imply
\begin{equation}\label{eq:p27-remainder-pre}
\begin{split}
&\|\Pi_n\partial_t^q\mathcal E(Q_1,Q_2)\Pi_{n'}\|_{2\to2}\le
CK_1K_2W_{q+b}
e^{-\rho\phi(|n-n'|)}
h(n-n')
\bigl(1-\widetilde\chi_j(n,n')\bigr).
\end{split}
\end{equation}
Write $d_*=|n-n'|$. On the remainder support,
\[
d_*>\frac{r_{nn'}}{10(j+2)},
\]
so
\[
X_{nn'}\le10(j+2)(1+d_*).
\]
Consequently, for every $N\in\mathbb N_0$,
\begin{equation}\label{eq:p27-tail-start}
X_{nn'}^{\tau N}e^{-\rho\phi(d_*)}
\le
[10(j+2)]^{\tau N}
(1+d_*)^{\tau N}e^{-\rho\phi(d_*)}.
\end{equation}

Recall that
\[
\phi(r)=[\log(R_\sigma+r)]^\sigma-(\log R_\sigma)^\sigma,
\qquad R_\sigma\ge1.
\]
Hence
\[
\phi(r)\ge[\log(1+r)]^\sigma-(\log R_\sigma)^\sigma.
\]
With $y=\log(1+d_*)\ge0$, this yields
\[
(1+d_*)^{\tau N}e^{-\rho\phi(d_*)}
\le C_{\sigma,\rho}e^{\tau Ny-\rho y^\sigma}.
\]
The moment estimate
\[
\sup_{y\ge0}\{\tau Ny-\rho y^\sigma\}
=
D_\sigma\rho^{\frac{-1}{\sigma-1}}\tau^pN^p,
\qquad
D_\sigma=(\sigma-1)\sigma^{-\frac{\sigma}{\sigma-1}},
\]
therefore gives
\begin{equation}\label{eq:p27-tail}
X_{nn'}^{\tau N}e^{-\rho\phi(d_*)}
\le
C_{\sigma,\rho}
\exp\!\left(
\tau N\log(10(j+2))+b_\rho N^p
\right),
\end{equation}
where
\[
b_\rho=D_\sigma\rho^{-\frac{1}{\sigma-1}}\tau^p.
\]

We then show that the fixed choices
\begin{equation}\label{eq:p27-kappa}
\kappa\ge2b_\rho
\end{equation}
and
\begin{equation}\label{eq:p27-M}
\log M\ge
\max\left(
0,\sup_{j\in\mathbb N_0}
\left\{
\tau\log(10(j+2))
-\frac{\kappa p}{2}j^{p-1}
\right\}
\right)
\end{equation}
imply
\begin{equation}\label{eq:p27-absorb}
W_{q+b}
\exp\!\left(
\tau N\log(10(j+2))+b_\rho N^p
\right)
\le W_{q+j+N}.
\end{equation}
The supremum in \eqref{eq:p27-M} is finite, since $p-1>0$
and $j^{p-1}$ grows faster than $\log(j+2)$.

Since $b\le j$ and $M\ge1$,
\begin{align*}
\frac{W_{q+j+N}}{W_{q+b}}
&=
M^{j+N-b}
\left(\frac{(q+j+N)!}{(q+b)!}\right)^2
e^{\kappa[(q+j+N)^p-(q+b)^p]}\ge
M^N e^{\kappa[(j+N)^p-j^p]}.
\end{align*}
For the last inequality, first replace $q+b$ by $q+j$,
and then use the monotonicity of
\[
x\longmapsto(x+N)^p-x^p
\qquad (x\ge0).
\]
Moreover, superadditivity and convexity give respectively
\[
(j+N)^p-j^p\ge N^p,
\qquad
(j+N)^p-j^p\ge pNj^{p-1}.
\]
Averaging these two inequalities yields
\[
(j+N)^p-j^p
\ge\frac12N^p+\frac p2Nj^{p-1}.
\]
Therefore
\[
\frac{W_{q+j+N}}{W_{q+b}}
\ge
\exp\!\left(
N\log M+\frac{\kappa}{2}N^p
+\frac{\kappa p}{2}Nj^{p-1}
\right).
\]
Conditions \eqref{eq:p27-kappa} and \eqref{eq:p27-M}
now give \eqref{eq:p27-absorb}.

Combining \eqref{eq:p27-remainder-pre},
\eqref{eq:p27-tail}, and \eqref{eq:p27-absorb}, we obtain,
for every $q,N\in\mathbb N_0$,
\[
\|\Pi_n\partial_t^q\mathcal E(Q_1,Q_2)\Pi_{n'}\|_{2\to2}
\le
CK_1K_2W_{q+j+N}
h(n-n')X_{nn'}^{-\tau N}.
\]
Thus we have
\[
|\mathcal E(Q_1,Q_2)|_j\le CK_1K_2.
\]
Together with the principal-term estimate , this proves
\[
\|\mathcal M(Q_1,Q_2)\|_{j,0}
+
|\mathcal E(Q_1,Q_2)|_j
\le
C\|Q_1\|_{j_1,0}\|Q_2\|_{j_2,0}.
\]
The parameters $M,\kappa$ are chosen once and for all,
independently of $q,N,j_1,j_2$.
\end{proof}

\begin{proposition}[Refined composition for mixed $\delta$]\label{2.8}
If $j_1\ge1,j_2\ge0$, $Q_1\in\mathcal{L}_\tau^{-j_1}(0)$, and $Q_2\in\mathcal{L}_\tau^{-j_2}(1)$, then
\[
 Q_1Q_2=\mathcal{M}(Q_1,Q_2)+\mathcal{E}(Q_1,Q_2),\quad
 Q_2Q_1=\mathcal{M}(Q_2,Q_1)+\mathcal{E}(Q_2,Q_1),
\]
both principal terms belong to $\widetilde{\mathcal{L}}_\tau^{-(j_1+j_2)}(0)$, both remainders belong to
$\mathcal{R}_{j_1+j_2}^{-\infty}$, and
\begin{equation}\label{eq:mixeddelta}
\begin{split}
&\|\mathcal M(Q_1,Q_2)\|_{j,0}
+\|\mathcal M(Q_2,Q_1)\|_{j,0}
+|\mathcal E(Q_1,Q_2)|_j
+|\mathcal E(Q_2,Q_1)|_j\le
C\|Q_1\|_{j_1,0}\|Q_2\|_{j_2,1}.
\end{split}
\end{equation}
\end{proposition}
\begin{proof}
Set
\[
j=j_1+j_2,\qquad
K_1=\|Q_1\|_{j_1,0},\qquad
K_2=\|Q_2\|_{j_2,1}.
\]
Since $j_1\ge1$ and $j_2\ge0$, we have $j\ge1$.
All constants below are independent of $j_1,j_2,q,N$, and $t$.

Leibniz's formula gives
\begin{equation}\label{eq:p28-product}
\begin{split}
\Pi_n\partial_t^q(Q_1Q_2)\Pi_{n'}
={}&
\sum_{\ell=0}^q\binom q\ell
\sum_{k\in\mathbb Z^d}
\bigl(\Pi_n\partial_t^\ell Q_1\Pi_k\bigr)
\bigl(\Pi_k\partial_t^{q-\ell}Q_2\Pi_{n'}\bigr).
\end{split}
\end{equation}
After dropping the input support indicators, the defining
matrix estimates imply
\begin{align*}
\|\Pi_n\partial_t^\ell Q_1\Pi_k\|_{2\to2}
&\le
K_1W_{\ell+j_1-1}
e^{-\rho\phi(|n-k|)}
h(n-k)X_{nk}^{-j_1\tau},\\
\|\Pi_k\partial_t^{q-\ell}Q_2\Pi_{n'}\|_{2\to2}
&\le
K_2W_{q-\ell+j_2}
e^{-\rho\phi(|k-n'|)}
h(k-n')X_{kn'}^{-j_2\tau}.
\end{align*}

On the support of the principal cutoff,
\[
|n-n'|\le\frac{r_{nn'}}{10(j+2)}.
\]
The radial estimate established in the proof of Proposition \ref{2.7}
therefore gives
\begin{equation}\label{eq:p28-radial}
X_{nk}^{-j_1\tau}X_{kn'}^{-j_2\tau}
\le e^{\tau/10}X_{nn'}^{-j\tau}.
\end{equation}
This also holds when $j_2=0$, in which case the second radial
factor equals one.\\
The defining estimates for the two factors involve the weights
\[
W_{\ell+j_1-1}
\qquad\text{and}\qquad
W_{q-\ell+j_2},
\]
respectively. The second expression remains valid when $j_2=0$,
since
\[
(j_2+1-1)_+=j_2.
\]
Applying the temporal convolution inequality (\ref{t1}) with
\[
b_1=j_1-1\ge0,\qquad b_2=j_2\ge0,
\]
we obtain
\begin{equation}\label{eq:p28-time}
\sum_{\ell=0}^q\binom q\ell
W_{\ell+j_1-1}W_{q-\ell+j_2}
\le 3W_{q+j-1}.
\end{equation}
Since $j\ge1$, the resulting index is
\[
q+j-1=q+(j-1)_+,
\]
which is precisely the temporal index required for the output
class $\mathcal L_\tau^{-j}(0)$.\\
Combining the spatial convolution estimate (\ref{s1}) with
\eqref{eq:p28-time}--\eqref{eq:p28-radial}, we obtain
\[
\begin{split}
&\|\Pi_n\partial_t^q\mathcal M(Q_1,Q_2)\Pi_{n'}\|_{2\to2}\le
CK_1K_2W_{q+j-1}
e^{-\rho\phi(|n-n'|)}
h(n-n')X_{nn'}^{-j\tau}
\widetilde\chi_j(n,n').
\end{split}
\]
Thus
\[
\mathcal M(Q_1,Q_2)
\in\widetilde{\mathcal L}_\tau^{-j}(0),
\qquad
\|\mathcal M(Q_1,Q_2)\|_{j,0}\le CK_1K_2.
\]

For the remainder, discard the radial factors using
\[
X_{nk}^{-j_1\tau}X_{kn'}^{-j_2\tau}\le1.
\]
The same temporal and spatial convolution estimates yield
\begin{equation}\label{eq:p28-remainder-pre}
\begin{split}
&\|\Pi_n\partial_t^q\mathcal E(Q_1,Q_2)\Pi_{n'}\|_{2\to2}\le
CK_1K_2W_{q+j-1}
e^{-\rho\phi(|n-n'|)}
h(n-n')
\bigl(1-\widetilde\chi_j(n,n')\bigr).
\end{split}
\end{equation}
On the remainder support, setting $d_*=|n-n'|$, we have
\[
X_{nn'}\le10(j+2)(1+d_*).
\]
Consequently, for every $N\in\mathbb N_0$,
\[
X_{nn'}^{\tau N}e^{-\rho\phi(d_*)}
\le
C_{\sigma,\rho}
\exp\!\left(
\tau N\log(10(j+2))+b_\rho N^p
\right),
\]
where
\[
b_\rho=D_\sigma\rho^{-\frac{1}{\sigma-1}}\tau^p.
\]

Use the same fixed parameter choices as in Proposition \ref{2.7}:
\[
\kappa\ge2b_\rho,
\]
and
\[
\log M\ge
\max\left(
0,\sup_{r\in\mathbb N_0}
\left\{
\tau\log(10(r+2))
-\frac{\kappa p}{2}r^{p-1}
\right\}
\right).
\]
The absorption estimate proved there states that, for
$0\le b\le j$,
\[
W_{q+b}
\exp\!\left(
\tau N\log(10(j+2))+b_\rho N^p
\right)
\le W_{q+j+N}.
\]
Taking $b=j-1$, which is nonnegative, and applying this inequality
to \eqref{eq:p28-remainder-pre}, we obtain
\[
\|\Pi_n\partial_t^q\mathcal E(Q_1,Q_2)\Pi_{n'}\|_{2\to2}
\le
CK_1K_2W_{q+j+N}
h(n-n')X_{nn'}^{-\tau N}.
\]
This holds for all $q,N\in\mathbb N_0$, uniformly in time.
Hence
\[
\mathcal E(Q_1,Q_2)\in\mathcal R_j^{-\infty},
\qquad
|\mathcal E(Q_1,Q_2)|_j\le CK_1K_2.
\]
Taking adjoints preserves the operator classes and their
seminorms. Therefore
\[
Q_1^*\in\mathcal L_\tau^{-j_1}(0),
\qquad
Q_2^*\in\mathcal L_\tau^{-j_2}(1).
\]
Apply the result just proved to $Q_1^*Q_2^*$ and take adjoints.
Since the output cutoff is symmetric in $n,n'$,
\[
\mathcal M(Q_1^*,Q_2^*)^*
=\mathcal M(Q_2,Q_1),
\]
and likewise
\[
\mathcal E(Q_1^*,Q_2^*)^*
=\mathcal E(Q_2,Q_1).
\]
Thus the reversed product has the same class memberships and
seminorm bounds.

Combining the two estimates and enlarging $C$ gives
\[
\begin{split}
&\|\mathcal M(Q_1,Q_2)\|_{j,0}
+\|\mathcal M(Q_2,Q_1)\|_{j,0}
+|\mathcal E(Q_1,Q_2)|_j
+|\mathcal E(Q_2,Q_1)|_j\le
C\|Q_1\|_{j_1,0}\|Q_2\|_{j_2,1}.
\end{split}
\]
All parameters remain fixed independently of the grades.
\end{proof}

\medskip
\noindent
\textbf{Remark}
The proof includes $j_2=0$. The hypothesis $j_1\ge1$ is used in
the identity
\[
(j_1-1)_++j_2=j_1+j_2-1.
\]
If $j_1=0$ and $j_2\ge1$, the left-hand side becomes $j_2$,
whereas the target class $(0)$ requires the index $j_2-1$.
The resulting extra temporal weight cannot be absorbed by a
constant uniform in the number of time derivatives. Thus this
argument does not justify removing the positive-order
hypothesis on the class $(0)$ factor.

\begin{corollary}[Symmetrized binary composition]\label{5}
Under the hypotheses of Proposition \ref{2.7}, there exist self-adjoint operators $\mathcal{M}_s\in\widetilde{\mathcal{L}}_\tau^{-j}(0)$ and
$\mathcal{E}_s\in\mathcal{R}_j^{-\infty}$ such that
\[
 Q_1Q_2+(Q_1Q_2)^*=\mathcal{M}_s+\mathcal{E}_s,
 \qquad \|\mathcal{M}_s\|_{j,0}+|\mathcal{E}_s|_j
 \le C\|Q_1\|_{j_1,0}\|Q_2\|_{j_2,0}.
\]
\end{corollary}
\begin{proof}
Take $\mathcal{M}_s=\mathcal{M}+\mathcal{M}^*$ and $\mathcal{E}_s=\mathcal{E}+\mathcal{E}^*$. Adjoint seminorms are equal, so no separate adjoint norms need to be listed.
\end{proof}

\begin{corollary}[Symmetrized mixed composition]\label{6}
Under the hypotheses of Proposition \ref{2.8}, both $Q_1Q_2+(Q_1Q_2)^*$ and
$Q_2Q_1+(Q_2Q_1)^*$ admit the preceding self-adjoint decomposition. The principal term belongs to
$\widetilde{\mathcal{L}}_\tau^{-j}(0)$, the remainder belongs to $\mathcal{R}_j^{-\infty}$, and the seminorms are bounded by
$C\|Q_1\|_{j_1,0}\|Q_2\|_{j_2,1}$.
\end{corollary}

\begin{proposition}[Narrowing the support cutoff]\label{2.9}
If $j\ge1$ and $Q\in\mathcal{L}_\tau^{-j}(1)$, define
\[
 \Pi_n\widetilde Q\Pi_{n'}=\Pi_nQ\Pi_{n'}\widetilde\chi_j(n,n'),\qquad
 \widetilde R=Q-\widetilde Q.
\]
Then
\[
 \widetilde Q\in\widetilde{\mathcal{L}}_\tau^{-j}(1),\quad
 \|\widetilde Q\|_{j,1}\le\|Q\|_{j,1},\quad
 \widetilde R\in\mathcal{R}_j^{-\infty},\quad
 |\widetilde R|_j\le C\|Q\|_{j,1}.
\]
If $Q$ is self-adjoint, both terms are self-adjoint.
\end{proposition}
\begin{proof}
The proof uses the same cutoff-tail estimate as Proposition \ref{2.7}.
\end{proof}

\begin{proposition}[The two-sided ideal property of smoothing remainders]\label{2.10}
Let $0<\tau,\tau'\le1$, $Q\in\mathcal{L}_{\tau'}^{-j_1}(0)$, and
$R\in\mathcal{R}_{j_2}^{-\infty}(M,\kappa,\tau)$. Then
\begin{equation}\label{eq:ideal}
 QR,RQ\in\mathcal{R}_{j_1+j_2}^{-\infty}(M,\kappa,\tau),\qquad
 |QR|_{j_1+j_2}+|RQ|_{j_1+j_2}
 \le C\|Q\|_{j_1,0}|R|_{j_2}.
\end{equation}
\end{proposition}

\begin{proof}
Consider $RQ$. For a fixed intermediate frequency $\ell$, first combine the time derivatives using~(\ref{t1}), and set
$T=q+j_1+j_2$. The matrix product is then bounded by a constant times
\[
 h(n-\ell)h(\ell-n')e^{-\rho\phi(|\ell-n'|)}G_T(X_{n\ell}).
\]
If $X_{nn'}\le X_{n\ell}$, use the fact that $G_T$ is nonincreasing. Otherwise, the near-diagonal support of $Q$ gives
\[
 \frac{X_{nn'}}{X_{n\ell}}\le\frac{10}9,\qquad
 X_{nn'}-X_{n\ell}\le|\ell-n'|.
\]
Apply~\eqref{eq:transfer}; the resulting growth is absorbed by the spatial exponential decay of $Q$. Summing over $\ell$ gives
$C h(n-n')G_T(X_{nn'})$. For every $N$,
$G_T(X_{nn'})\le W_{T+N}X_{nn'}^{-\tau N}$, which is the required remainder bound. The proof for $QR$ is identical.
\end{proof}

\begin{proposition}[Ternary composition]\label{2.11}
If $Q_i\in\mathcal{L}_\tau^{-j_i}(0)$, $j_i\ge0$, and $j=j_1+j_2+j_3$, then
\[
 Q_1Q_2Q_3=Q+R,\quad
 Q\in\widetilde{\mathcal{L}}_\tau^{-j}(0),\quad R\in\mathcal{R}_j^{-\infty},
\]
and
\[
 \|Q\|_{j,0}+|R|_j\le C\prod_{i=1}^3\|Q_i\|_{j_i,0}.
\]
\end{proposition}
\begin{proof}
Apply Proposition \ref{2.7} to the first two factors, then Proposition \ref{2.7} to the principal term and Proposition \ref{2.10} to the remainder.
\end{proof}
\begin{proposition}[Placing the potential in the operator classes]\label{3.11}
Under (\ref{eq:DC}), choose $0<\rho<\rho_0$ and the fixed parameters $M,\kappa$. The multiplication operator can be written as
\begin{equation}\label{eq:Vsplit}
 V=Q_V+R_V,
\end{equation}
where
\begin{equation}\label{eq:Vsplit1}
 Q_V\in\widetilde{\mathcal{L}}_\tau^0(0),\quad R_V\in\mathcal{R}_0^{-\infty},
\end{equation}
and
\begin{equation}\label{eq:Vsplit2}
\|Q_V\|_{0,0}+|R_V|_0\le H_{0}.
\end{equation}
If $V$ is real-valued, then $Q_V $ and $ R_V$ are self-adjoint. Here $H_{0}$ depends only on the fixed data.
\end{proposition}

\begin{proof}
We write the spatial decay assumption in the form
\begin{equation}\label{eq:p317-assumption}
\sup_{t\in\mathbb R}
|\partial_t^q\widehat V(k,t)|
\le
A e^{Bq^p}e^{-\rho_0\phi(|k|)},
\qquad
p=\frac{\sigma}{\sigma-1}>1,
\end{equation}
where
\[
\phi(r)
=
[\log(R_\sigma+r)]^\sigma-(\log R_\sigma)^\sigma,
\qquad
R_\sigma=\max(e,e^{\sigma-1}).
\]
Fix $0<\rho<\rho_0$, and recall
\[
r_{nn'}=\max(|n|,|n'|),\qquad
X_{nn'}=1+r_{nn'},\qquad
h(k)=\langle k\rangle^{-(d+2)}.
\]

Since $V(t)$ acts by multiplication,
\[
\Pi_nV(t)\Pi_{n'}
=
\widehat V(n-n',t)\,
\langle\,\cdot\,,e_{n'}\rangle e_n.
\]
Consequently,
\[
\|\Pi_n\partial_t^qV(t)\Pi_{n'}\|_{2\to2}
=
|\partial_t^q\widehat V(n-n',t)|.
\]
Define the symmetric cutoff
\[
\chi_0(n,n')
=
\mathbf 1_{\{|n-n'|\le r_{nn'}/20\}},
\]
and set
\begin{align*}
\Pi_nQ_V(t)\Pi_{n'}
&=\chi_0(n,n')\Pi_nV(t)\Pi_{n'},\\
\Pi_nR_V(t)\Pi_{n'}
&=(1-\chi_0(n,n'))\Pi_nV(t)\Pi_{n'}.
\end{align*}
Thus $V=Q_V+R_V$ at the level of Fourier matrices.

Since $\sigma>1$,
\[
\frac{\phi(r)}{\log\langle r\rangle}\longrightarrow\infty
\qquad\text{as }r\to\infty.
\]
Hence
\[
C_{\mathrm{sp}}
:=
\sup_{k\in\mathbb Z^d}
\langle k\rangle^{d+2}
e^{-(\rho_0-\rho)\phi(|k|)}
<\infty.
\]
It follows that
\[
e^{-\rho_0\phi(|k|)}
\le
C_{\mathrm{sp}}h(k)e^{-\rho\phi(|k|)}.
\]
Therefore \eqref{eq:p317-assumption} implies
\begin{equation}\label{eq:p317-matrix}
\begin{split}
\sup_t
\|\Pi_n\partial_t^qV(t)\Pi_{n'}\|_{2\to2}
\le{}&
AC_{\mathrm{sp}}e^{Bq^p}
h(n-n')e^{-\rho\phi(|n-n'|)}.
\end{split}
\end{equation}

Choose $M\ge1$ and $\kappa\ge B$. Then
\[
e^{Bq^p}\le W_q,
\qquad
W_b=M^b(b!)^2e^{\kappa b^p}.
\]
Using \eqref{eq:p317-matrix}, we obtain
\[
\begin{split}
\sup_t
\|\Pi_n\partial_t^qQ_V(t)\Pi_{n'}\|_{2\to2}
\le{}&
AC_{\mathrm{sp}}W_q
e^{-\rho\phi(|n-n'|)}h(n-n')\times
\mathbf 1_{\{|n-n'|\le r_{nn'}/20\}}.
\end{split}
\]
At grade zero, the narrower cutoff is exactly
\[
\frac{r_{nn'}}{10(0+2)}=\frac{r_{nn'}}{20},
\]
and the temporal weight is
\[
W_{q+(0+0-1)_+}=W_q.
\]
Thus
\begin{equation}\label{eq:p317-principal}
Q_V\in\widetilde{\mathcal L}_\tau^0(0),
\qquad
\|Q_V\|_{0,0}\le AC_{\mathrm{sp}}.
\end{equation}

On the support of $R_V$, we have
\[
|n-n'|>\frac{r_{nn'}}{20}.
\]
Writing $d_*=|n-n'|$, we obtain
\[
X_{nn'}\le20(1+d_*).
\]
Consequently, for every $N\in\mathbb N_0$,
\[
X_{nn'}^{\tau N}e^{-\rho\phi(d_*)}
\le
20^{\tau N}(1+d_*)^{\tau N}e^{-\rho\phi(d_*)}.
\]
Since $R_\sigma\ge1$,
\[
\phi(r)
\ge
[\log(1+r)]^\sigma-(\log R_\sigma)^\sigma.
\]
Setting $y=\log(1+d_*)\ge0$ therefore gives
\[
(1+d_*)^{\tau N}e^{-\rho\phi(d_*)}
\le
e^{\rho(\log R_\sigma)^\sigma}
e^{\tau Ny-\rho y^\sigma}.
\]
The moment estimate yields
\[
\sup_{y\ge0}\{\tau Ny-\rho y^\sigma\}
=
D_\sigma\rho^{-1/(\sigma-1)}\tau^pN^p,
\qquad
D_\sigma=(\sigma-1)\sigma^{-\sigma/(\sigma-1)}.
\]
Put
\[
b_\rho=D_\sigma\rho^{-1/(\sigma-1)}\tau^p.
\]
We conclude that
\begin{equation}\label{eq:p317-tail}
X_{nn'}^{\tau N}e^{-\rho\phi(d_*)}
\le
C_{\sigma,\rho}
e^{\tau N\log20+b_\rho N^p}.
\end{equation}

Combining \eqref{eq:p317-matrix} and \eqref{eq:p317-tail}
on the remainder support gives
\[
\begin{split}
\sup_t
\|\Pi_n\partial_t^qR_V(t)\Pi_{n'}\|_{2\to2}
\le{}&
CA\,h(n-n')X_{nn'}^{-\tau N}\times
e^{Bq^p+b_\rho N^p+\tau N\log20}.
\end{split}
\]
Choose the parameters once and for all so that
\[
\kappa\ge\max(B,b_\rho),
\qquad
\log M\ge\tau\log20.
\]
Since $p>1$,
\[
q^p+N^p\le(q+N)^p.
\]
Hence
\[
Bq^p+b_\rho N^p+\tau N\log20
\le
\kappa(q+N)^p+(q+N)\log M,
\]
and therefore
\[
e^{Bq^p+b_\rho N^p+\tau N\log20}
\le W_{q+N}.
\]
It follows that
\begin{equation}\label{eq:p317-remainder}
\sup_t
\|\Pi_n\partial_t^qR_V(t)\Pi_{n'}\|_{2\to2}
\le
CAW_{q+N}h(n-n')X_{nn'}^{-\tau N}
\end{equation}
for all $q,N\in\mathbb N_0$. This is precisely the defining
estimate for $\mathcal R_0^{-\infty}$, so
\[
R_V\in\mathcal R_0^{-\infty},
\qquad
|R_V|_0\le CA.
\]

Since $V$ is real-valued,
\[
\widehat V(-k,t)=\overline{\widehat V(k,t)}.
\]
Moreover,
\[
\chi_0(n,n')=\chi_0(n',n).
\]
Thus
\[
(Q_V)_{n'n}^*=(Q_V)_{nn'},
\qquad
(R_V)_{n'n}^*=(R_V)_{nn'}.
\]
The preceding estimates and $h\in\ell^1(\mathbb Z^d)$ show
that both operators are bounded on $L^2$. Their matrix symmetry
therefore implies
\[
Q_V^*=Q_V,\qquad R_V^*=R_V.
\]
Combining \eqref{eq:p317-principal} and
\eqref{eq:p317-remainder}, and absorbing the fixed constants
into $H_{0}$, gives
\[
\|Q_V\|_{0,0}+|R_V|_0\le H_{0}.
\]
This completes the proof.
\end{proof}

First choose $0<\tau\le\tau_0$ and $0<\rho<\rho_0$. Then choose $\kappa$ once and for all, larger than
the temporal constant $B$, the threshold required in~\eqref{eq:threshold}, and $b_g$ from Lemma \ref{D}. Finally choose $M$ sufficiently large. These choices are independent of the normalization order.

\begin{lemma}[Graded expansion under logarithmic Fourier decay]\label{L3.1*}Set
\[
P_0=i\partial_t-\Delta.
\]
Suppose that the real-valued multiplication potential admits a
decomposition
\[
V=Q_V+R_V,
\qquad
Q_V\in\widetilde{\mathcal L}_\tau^0(0),
\qquad
R_V\in\mathcal R_0^{-\infty},
\]
where $Q_V$ and $R_V$ are self-adjoint. Choose
\[
H\ge1+\|Q_V\|_{0,0}+|R_V|_0.
\]
Let $K$ be a large constant. Suppose that there are two sequences $Q'_r,Q''_r$, $1\le r\le m$, satisfying
\[
Q'_r\in\mathcal L_\tau^{-r}(0),
(Q'_r)^*=-Q'_r,  [\Delta,Q'_r]\in
\mathcal L_\tau^{-(r-1)}(1),
\]
\[
Q''_r\in\mathcal L_\tau^{-(r+1)}(0),
(Q''_r)^*=Q''_r, [\Delta,Q''_r]\in
\mathcal L_\tau^{-r}(0),
\]
and the quantitative estimates
\begin{align}
\|Q'_r\|_{r,0}
&\le
\frac{K^{r-\frac12}H^r}{r^2},
&
\|[\Delta,Q'_r]\|_{r-1,1}
&\le
\frac{K^{r-1}H^r}{r^2},
\label{eq:log-graded-prime}
\\
\|Q''_r\|_{r+1,0}
&\le
\frac{K^{r+\frac12}H^{r+1}}{(r+1)^2},
&
\|[\Delta,Q''_r]\|_{r,0}
&\le
\frac{K^rH^{r+1}}{(r+1)^2}.
\label{eq:log-graded-double-prime}
\end{align}

Define
\[
\mathcal Q_m=\sum_{r=1}^m(Q'_r+Q''_r).
\]
Then there exist self-adjoint operators
\[
S_j\in
\widetilde{\mathcal L}_\tau^{-(j+1)}(0),
\qquad
R_j\in\mathcal R^{-\infty}_{j+1},
\qquad 1\le j\le 2m+1,
\]
\[
\widetilde S_j\in
\widetilde{\mathcal L}_\tau^{-(j-1)}(1),
\qquad 2\le j\le m+1,
\]
and
\[
T_j\in
\widetilde{\mathcal L}_\tau^{-(j-1)}(0),
\qquad
U_j\in\mathcal R^{-\infty}_{j-1},
\qquad 2\le j\le 2m+3,
\]
such that
\begin{align}
&(I+\mathcal Q_m)^*(P_0+V)(I+\mathcal Q_m)
\nonumber\\
&\quad=
P_0+V
+\sum_{r=1}^m[Q'_r,\Delta]
+\sum_{r=1}^m(Q''_rP_0+P_0Q''_r)
\nonumber\\
&\qquad\quad
+\frac12\sum_{j=1}^{2m+1}(S_jP_0+P_0S_j)
+\frac12\sum_{j=1}^{2m+1}(R_jP_0+P_0R_j)
\nonumber\\
&\qquad\quad
+\sum_{j=2}^{m+1}\widetilde S_j
+\sum_{j=2}^{2m+3}(T_j+U_j).
\label{eq:log-graded-expansion}
\end{align}
This identity initially holds on smooth functions compactly supported
in time.

Moreover,
\[
[\Delta,S_j]\in
\widetilde{\mathcal L}_\tau^{-j}(0),
\]
and
\begin{align}
\|S_j\|_{j+1,0}+|R_j|_{j+1}
&\le
C\frac{K^jH^{j+1}}{(j+1)^2},
\label{eq:log-S-estimate}
\\
\|[\Delta,S_j]\|_{j,0}
&\le
C\frac{K^{j-\frac12}H^{j+1}}{(j+1)^2},
\label{eq:log-comm-S-estimate}
\\
\|\widetilde S_j\|_{j-1,1}
&\le
C\frac{K^{j-\frac32}H^{j-1}}{j^2},
\label{eq:log-tilde-S-estimate}
\\
\|T_j\|_{j-1,0}+|U_j|_{j-1}
&\le
C\frac{K^{j-\frac32}H^j}{j^2}.
\label{eq:log-T-U-estimate}
\end{align}
The constant $C$ depends only on the fixed parameters of the operator
calculus, and is independent of $m,j,K,H$.
\end{lemma}

\begin{proof}
We use the ordinary composition theorem, the mixed composition theorem,
the cutoff-narrowing property, and the remainder ideal property.
In particular, the mixed composition theorem applies to
\[
\mathcal L_\tau^{-u}(0)\,
\mathcal L_\tau^{-v}(1),
\qquad u\ge1,\quad v\ge0,
\]
in either order. Its principal part belongs to $\widetilde{\mathcal L}_\tau^{-(u+v)}(0) $, and its remainder belongs to $ \mathcal R^{-\infty}_{u+v}$.
The endpoint $v=0$ is included.

All sums over $r,s$ below are restricted to $1\le r,s\le m$.
An empty sum is understood to be zero.

\medskip
\noindent
\textbf{Step 1. The exact algebraic identity.}

Since
\[
(Q'_r)^*=-Q'_r,
\qquad
(Q''_r)^*=Q''_r,
\]
the terms linear in $\mathcal Q_m$ satisfy
\begin{align*}
\mathcal Q_m^*P_0+P_0\mathcal Q_m
&=
\sum_{r=1}^m[Q'_r,\Delta]
-\sum_{r=1}^m[Q'_r,i\partial_t]
\\
&\quad+
\sum_{r=1}^m(Q''_rP_0+P_0Q''_r).
\end{align*}
For the quadratic term, we use
\begin{align*}
\mathcal Q_m^*P_0\mathcal Q_m
&=
\frac12\bigl(
\mathcal Q_m^*\mathcal Q_mP_0
+P_0\mathcal Q_m^*\mathcal Q_m
\bigr)
\\
&\quad+
\frac12\bigl(
\mathcal Q_m^*[P_0,\mathcal Q_m]
+[\mathcal Q_m^*,P_0]\mathcal Q_m
\bigr).
\end{align*}
Consequently,
\begin{align}
&(I+\mathcal Q_m)^*(P_0+V)(I+\mathcal Q_m)
\nonumber\\
&=
P_0+V
+\sum_{r=1}^m[Q'_r,\Delta]
+\sum_{r=1}^m(Q''_rP_0+P_0Q''_r)
-\sum_{r=1}^m[Q'_r,i\partial_t]
\nonumber\\
&\quad+
\frac12\bigl(
\mathcal Q_m^*\mathcal Q_mP_0
+P_0\mathcal Q_m^*\mathcal Q_m
\bigr)
\nonumber\\
&\quad+
\frac12\bigl(
\mathcal Q_m^*[i\partial_t,\mathcal Q_m]
+[\mathcal Q_m^*,i\partial_t]\mathcal Q_m
\bigr)
\nonumber\\
&\quad+
\frac12\bigl(
\mathcal Q_m^*[-\Delta,\mathcal Q_m]
+[\mathcal Q_m^*,-\Delta]\mathcal Q_m
\bigr)
\nonumber\\
&\quad+
\mathcal Q_m^*V+V\mathcal Q_m
+\mathcal Q_m^*V\mathcal Q_m.
\label{eq:log-exact-identity}
\end{align}

\medskip
\noindent
\textbf{Step 2. Decomposition of
$\mathcal Q_m^*\mathcal Q_m$.}

We group the products according to their total smoothing order:
\begin{align}
\mathcal Q_m^*\mathcal Q_m
&=
\sum_{j=1}^{2m-1}
\sum_{r+s=j+1}(Q'_r)^*Q'_s
\nonumber\\
&\quad+
\sum_{j=2}^{2m}
\sum_{r+s=j}
\bigl((Q'_r)^*Q''_s+(Q''_s)^*Q'_r\bigr)
\nonumber\\
&\quad+
\sum_{j=3}^{2m+1}
\sum_{r+s=j-1}(Q''_r)^*Q''_s.
\label{eq:log-QstarQ-grouping}
\end{align}
For each fixed $j$, every product in the corresponding group has
total smoothing order $(j+1)\tau$.

Define $S_j$ by retaining, in this group, the matrix entries satisfying
\[
|n-n'|
\le
\frac{\max(|n|,|n'|)}{10(j+3)}.
\]
Let $R_j$ be the complementary part.
The ordinary composition theorem gives
\[
S_j\in
\widetilde{\mathcal L}_\tau^{-(j+1)}(0),
\qquad
R_j\in\mathcal R^{-\infty}_{j+1},
\]
and
\[
\mathcal Q_m^*\mathcal Q_m
=
\sum_{j=1}^{2m+1}(S_j+R_j).
\]

The three kinds of products in
\eqref{eq:log-QstarQ-grouping} have respective majorants
\[
\frac{K^jH^{j+1}}{r^2s^2},
\qquad
\frac{K^jH^{j+1}}{r^2(s+1)^2},
\qquad
\frac{K^jH^{j+1}}{(r+1)^2(s+1)^2}.
\]
We repeatedly use
\[
\sum_{\ell=1}^{N-1}
\frac{1}{\ell^2(N-\ell)^2}
\le \frac{C}{N^2},
\qquad N\ge2.
\]
Indeed, on $\ell\le N/2$ one has
$(N-\ell)^{-2}\le4N^{-2}$, and the other half of the sum
is treated symmetrically.
The shifted versions follow by changing the summation variable.
Thus
\[
\|S_j\|_{j+1,0}+|R_j|_{j+1}
\le
C\frac{K^jH^{j+1}}{(j+1)^2}.
\]

Each group in \eqref{eq:log-QstarQ-grouping} is self-adjoint.
Since the cutoff is symmetric in $n,n'$, both $S_j$ and $R_j$
are self-adjoint.

\medskip
\noindent
\textbf{Step 3. The commutator estimate for $S_j$.}

For the first type of product,
\[
[\Delta,(Q'_r)^*Q'_s]
=
[\Delta,(Q'_r)^*]Q'_s
+(Q'_r)^*[\Delta,Q'_s].
\]
Here
\[
[\Delta,(Q'_r)^*]
\in\mathcal L_\tau^{-(r-1)}(1),
\qquad
Q'_s\in\mathcal L_\tau^{-s}(0),
\]
and similarly for the second term.
Since $r,s\ge1$, the mixed composition theorem applies to both
products, including $r=1$ or $s=1$.
Their principal parts have class
\[
\widetilde{\mathcal L}_\tau^{-(r+s-1)}(0).
\]
Under $r+s=j+1$, this is
$\widetilde{\mathcal L}_\tau^{-j}(0)$.

For the second type,
\[
[\Delta,(Q'_r)^*Q''_s]
=
[\Delta,(Q'_r)^*]Q''_s
+(Q'_r)^*[\Delta,Q''_s].
\]
The first product is handled by mixed composition, and the second
by ordinary composition. Both have total smoothing order
\[
(r+s)\tau=j\tau.
\]

For the third type,
\[
[\Delta,(Q''_r)^*Q''_s]
=
[\Delta,(Q''_r)^*]Q''_s
+(Q''_r)^*[\Delta,Q''_s].
\]
Ordinary composition gives total smoothing order
\[
(r+s+1)\tau=j\tau.
\]

It remains to check that these estimates apply to the actual
commutator of $S_j$.
The Laplacian is diagonal in the Fourier basis, so taking a commutator
with $\Delta$ commutes with the matrix cutoff defining $S_j$.
For total order $j$, the composition decomposition uses the cutoff
\[
|n-n'|
\le
\frac{\max(|n|,|n'|)}{10(j+2)}.
\]
The cutoff defining $S_j$ is narrower:
\[
|n-n'|
\le
\frac{\max(|n|,|n'|)}{10(j+3)}.
\]
Therefore, after applying the latter cutoff, the complementary
remainders from these product decompositions vanish.
This proves
\[
[\Delta,S_j]\in
\widetilde{\mathcal L}_\tau^{-j}(0).
\]

In each differentiated product, replacing one factor by its
commutator decreases the power of $K$ by $1/2$ and leaves the
power of $H$ unchanged.
The same discrete convolution estimate as above yields
\[
\|[\Delta,S_j]\|_{j,0}
\le
C\frac{K^{j-\frac12}H^{j+1}}{(j+1)^2}.
\]

\medskip
\noindent
\textbf{Step 4. The linear time-derivative terms.}

For $2\le j\le m+1$, by Proposition \ref{2.4} and \ref{2.9}, one has
\[
-[Q'_{j-1},i\partial_t]
=
i\partial_tQ'_{j-1}
\in
\mathcal L_\tau^{-(j-1)}(1).
\]
By cutoff narrowing, write
\[
i\partial_tQ'_{j-1}
=
\widetilde S_j+E_j^{\mathrm{time}},
\]
where
\[
\widetilde S_j\in
\widetilde{\mathcal L}_\tau^{-(j-1)}(1),
\qquad
E_j^{\mathrm{time}}\in
\mathcal R^{-\infty}_{j-1}.
\]
The time-derivative estimate and
\eqref{eq:log-graded-prime} give
\[
\|\widetilde S_j\|_{j-1,1}
+|E_j^{\mathrm{time}}|_{j-1}
\le
C\frac{K^{j-\frac32}H^{j-1}}{(j-1)^2}
\le
C\frac{K^{j-\frac32}H^{j-1}}{j^2}.
\]
Here $j/(j-1)\le2$.

Since $Q'_{j-1}$ is skew-adjoint,
$i\partial_tQ'_{j-1}$ is self-adjoint.
The symmetric cutoff therefore makes both pieces self-adjoint.
The remainders $E_j^{\mathrm{time}}$ will be included in $U_j$.

\medskip
\noindent
\textbf{Step 5. The quadratic time-derivative terms.}

Consider
\[
\frac12\bigl(
\mathcal Q_m^*[i\partial_t,\mathcal Q_m]
+[\mathcal Q_m^*,i\partial_t]\mathcal Q_m
\bigr).
\]
The second summand is the adjoint of the first.
The products in the first summand are grouped as follows:
\[
\begin{array}{c|c}
\text{product} & \text{index relation}\\ \hline
(Q'_r)^*[i\partial_t,Q'_s] & r+s=j-1\\
(Q'_r)^*[i\partial_t,Q''_s] & r+s=j-2\\
(Q''_r)^*[i\partial_t,Q'_s] & r+s=j-2\\
(Q''_r)^*[i\partial_t,Q''_s] & r+s=j-3.
\end{array}
\]
A time derivative changes the flag from $(0)$ to $(1)$ without
changing the smoothing order.
The undifferentiated factor always has positive smoothing order.
Hence mixed composition gives, in every row, a principal part in $ \widetilde{\mathcal L}_\tau^{-(j-1)}(0) $
and a remainder in $ \mathcal R_{j-1}^{-\infty}.$

For all four rows, the product of the assumed majorants has
numerator $ K^{j-2}H^{j-1}$. \\
Summing the corresponding denominators gives $ C\frac{K^{j-2}H^{j-1}}{j^2} $.\ Since $K,H\ge1$, this is bounded by
$
C\frac{K^{j-\frac32}H^j}{j^2}.
$\\
Taking each principal part together with its adjoint preserves
self-adjointness, and likewise for the remainders.

\medskip
\noindent
\textbf{Step 6. The quadratic Laplacian-commutator terms.}

Next consider
\[
\frac12\bigl(
\mathcal Q_m^*[-\Delta,\mathcal Q_m]
+[\mathcal Q_m^*,-\Delta]\mathcal Q_m
\bigr).
\]
Again the two summands are adjoints of one another.
The relevant grouping is
\[
\begin{array}{c|c}
\text{product} & \text{index relation}\\ \hline
(Q'_r)^*[-\Delta,Q'_s] & r+s=j\\
(Q'_r)^*[-\Delta,Q''_s] & r+s=j-1\\
(Q''_r)^*[-\Delta,Q'_s] & r+s=j-1\\
(Q''_r)^*[-\Delta,Q''_s] & r+s=j-2.
\end{array}
\]
For a commutator involving $Q'_s$, use $ [\Delta,Q'_s]\in \mathcal L_\tau^{-(s-1)}(1) $. Its other factor has positive smoothing order, so mixed composition
applies. This includes $s=1$, when the commutator has order zero.\\
For a commutator involving $Q''_s$, use
\[
[\Delta,Q''_s]\in
\mathcal L_\tau^{-s}(0),
\]
and apply ordinary composition.\ Every row thus produces a principal part in
$ \widetilde{\mathcal L}_\tau^{-(j-1)}(0)$ and a remainder in
$ \mathcal R_{j-1}^{-\infty}.$\\
In every case, the majorant numerator is $ K^{j-\frac32}H^j $. Summing over $r,s$ proves the bound
$ C\frac{K^{j-\frac32}H^j}{j^2}. $

\medskip
\noindent
\textbf{Step 7. The terms linear in the potential.}

Split
\[
\mathcal Q_m^*V+V\mathcal Q_m
=
\mathcal Q_m^*Q_V+Q_V\mathcal Q_m
+\mathcal Q_m^*R_V+R_V\mathcal Q_m.
\]
For the products containing $Q'_r$, assign $ j=r+1;$ for those containing $Q''_r$, assign $ j=r+2.$ Since $Q_V$ has order zero and flag $(0)$, ordinary composition
gives principal parts in $ \widetilde{\mathcal L}_\tau^{-(j-1)}(0) $ and remainders in
$ \mathcal R^{-\infty}_{j-1}.$\\
Follow from
\[
\|Q_V\|_{0,0}\le H,
\]
their bounds for the prime terms,
\[
\frac{K^{r-\frac12}H^{r+1}}{r^2}
\le
C\frac{K^{j-\frac32}H^j}{j^2},
\]
and for the double-prime terms,
\[
\frac{K^{r+\frac12}H^{r+2}}{(r+1)^2}
\le
C\frac{K^{j-\frac32}H^j}{j^2}.
\]
The remainder ideal property gives the same order and bound for
the products containing $R_V$.

\medskip
\noindent
\textbf{Step 8. The terms quadratic in the normalizing operators
and linear in the potential.}

Write
\[
\mathcal Q_m^*V\mathcal Q_m
=
\mathcal Q_m^*Q_V\mathcal Q_m
+\mathcal Q_m^*R_V\mathcal Q_m.
\]
For the first expression, use
\[
\begin{array}{c|c}
\text{product} & \text{index relation}\\ \hline
(Q'_r)^*Q_VQ'_s & r+s=j-1\\
(Q'_r)^*Q_VQ''_s & r+s=j-2\\
(Q''_r)^*Q_VQ'_s & r+s=j-2\\
(Q''_r)^*Q_VQ''_s & r+s=j-3.
\end{array}
\]
The triple-composition property gives principal parts in
$
\widetilde{\mathcal L}_\tau^{-(j-1)}(0)
$
and remainders in
$
\mathcal R_{j-1}^{-\infty}.
$\\
Each row has majorant numerator
$
K^{j-2}H^j.
$ Thus their total contribution is bounded by
\[
C\frac{K^{j-2}H^j}{j^2}
\le
C\frac{K^{j-\frac32}H^j}{j^2}.
\]
The terms containing $R_V$ belong entirely to
$\mathcal R_{j-1}^{-\infty}$, with the same estimate, by the
remainder ideal property.

The largest possible index occurs in the double-prime,
double-prime product with $r=s=m$, namely
\[
j=2m+3.
\]

\medskip
\noindent
\textbf{Step 9. Collection of terms.}

For each $2\le j\le2m+3$, let $T_j$ be the sum of all principal
parts of order $-(j-1)\tau$ constructed in Steps 5--8.
Let $U_j$ be the sum of the corresponding remainders, together
with $E_j^{\mathrm{time}}$ from Step 4 when $j\le m+1$.

There are only finitely many types of contributions, independently
of $m$ and $j$.
The estimates obtained above therefore imply
\[
\|T_j\|_{j-1,0}+|U_j|_{j-1}
\le
C\frac{K^{j-\frac32}H^j}{j^2}.
\]
The adjoint pairing and symmetric cutoff ensure that $T_j$ and
$U_j$ are self-adjoint.

Substituting these decompositions and
\[
\mathcal Q_m^*\mathcal Q_m
=
\sum_{j=1}^{2m+1}(S_j+R_j)
\]
into \eqref{eq:log-exact-identity} proves
\eqref{eq:log-graded-expansion} and all the asserted estimates.
\end{proof}

\begin{proposition}\label{2.14}Let $V(x,t)$ be real-valued, with
\[
V(x,t)=\sum_{k\in\mathbb Z^d}v_k(t)e^{ik\cdot x},
\]
and assume that, for some $\sigma>1$,
\[
\sup_{t\in\mathbb R}|\partial_t^qv_k(t)|
\le
A\exp\!\left(Bq^{\frac{\sigma}{\sigma-1}}\right)
\exp\!\left(-\rho_0[\log(R_{\sigma}+|k|)^\sigma-(\log R_{\sigma})^{\sigma}]\right),
\]
where $ q\ge0, R_{\sigma}=\max\{e,e^{\sigma}\}$.\ Suppose that the real-valued multiplication potential admits a
decomposition
\[
V=Q_V+R_V,
\qquad
Q_V\in\widetilde{\mathcal L}_\tau^0(0),
\qquad
R_V\in\mathcal R_0^{-\infty},
\]
where $Q_V$ and $R_V$ are self-adjoint. Choose
\[
H\ge1+\|Q_V\|_{0,0}+|R_V|_0.
\]
There exists $K\ge1$, independent of $m$, such that, for every
$m\ge1$, one can construct operators $Q'_r,Q''_r$, $1\le r\le m$,
satisfying
\[
Q'_r\in\mathcal L_\tau^{-r}(0),
\qquad
(Q'_r)^*=-Q'_r,
\qquad
[\Delta,Q'_r]\in\mathcal L_\tau^{-(r-1)}(1),
\]
and
\[
Q''_r\in\mathcal L_\tau^{-(r+1)}(0),
\qquad
(Q''_r)^*=Q''_r,
\qquad
[\Delta,Q''_r]\in\mathcal L_\tau^{-r}(0).
\]
Moreover,
\begin{align}
\|Q'_r\|_{r,0}
&\le\frac{K^{r-\frac12}H^r}{r^2},
&
\|[\Delta,Q'_r]\|_{r-1,1}
&\le\frac{K^{r-1}H^r}{r^2},
\label{eq:normalization-prime-bounds}
\\
\|Q''_r\|_{r+1,0}
&\le\frac{K^{r+\frac12}H^{r+1}}{(r+1)^2},
&
\|[\Delta,Q''_r]\|_{r,0}
&\le\frac{K^rH^{r+1}}{(r+1)^2}.
\label{eq:normalization-double-bounds}
\end{align}

Set
\[
P_0=i\partial_t-\Delta,
\qquad
\mathcal Q_m=\sum_{r=1}^m(Q'_r+Q''_r).
\]
Then
\begin{align}
&(I+\mathcal Q_m)^*(P_0+V)(I+\mathcal Q_m)
\label{2.14*}\\
&=
P_0+V^{(m)}
+\frac12\sum_{j=m+1}^{2m+1}(S_jP_0+P_0S_j)
\nonumber\\
&\quad+
\frac12\sum_{j=1}^{2m+1}(R_jP_0+P_0R_j)
+\widetilde S_{m+1}
\nonumber\\
&\quad+
\sum_{j=m+1}^{2m+3}T_j
+\sum_{j=2}^{2m+3}U_j
+\sum_{j=1}^{m}\widehat R_j.
\label{eq:normalization-final}
\end{align}
The identity initially holds on smooth functions compactly
supported in time.

All coefficients on the right-hand side are self-adjoint.
The operator $V^{(m)}$ is block diagonal and satisfies
\[
[V^{(m)},\widetilde\Delta]=0.
\]
For each fixed $m$, it is uniformly bounded in time on $L^2$.

The remaining coefficients satisfy
\begin{align*}
&S_j\in\widetilde{\mathcal L}_\tau^{-(j+1)}(0),
\qquad
[\Delta,S_j]\in\widetilde{\mathcal L}_\tau^{-j}(0),
&&m+1\le j\le2m+1,\\
&R_j\in\mathcal R_{j+1}^{-\infty},
&&1\le j\le2m+1,\\
&\widetilde S_{m+1}\in
\widetilde{\mathcal L}_\tau^{-m}(1),\\
&T_j\in\widetilde{\mathcal L}_\tau^{-(j-1)}(0),
&&m+1\le j\le2m+3,\\
&U_j\in\mathcal R_{j-1}^{-\infty},
&&2\le j\le2m+3,\\
&\widehat R_j\in\mathcal R_{j-1}^{-\infty},
&&1\le j\le m.
\end{align*}
More precisely,
\begin{align}
\|S_j\|_{j+1,0}+|R_j|_{j+1}
&\le C\frac{K^jH^{j+1}}{(j+1)^2},
\label{eq:normalization-SR-bounds}\\
\|[\Delta,S_j]\|_{j,0}
&\le C\frac{K^{j-\frac12}H^{j+1}}{(j+1)^2},
\label{eq:normalization-commS-bound}\\
\|\widetilde S_{m+1}\|_{m,1}
&\le C\frac{K^{m-\frac12}H^m}{(m+1)^2},
\label{eq:normalization-time-bound}\\
\|T_j\|_{j-1,0}+|U_j|_{j-1}
&\le C\frac{K^{j-\frac32}H^j}{j^2},
\label{eq:normalization-TU-bounds}\\
|\widehat R_1|_0&\le CH,
\label{eq:normalization-Rhat-first}\\
|\widehat R_j|_{j-1}
&\le C\frac{K^{j-\frac32}H^j}{j^2},
\qquad 2\le j\le m.
\label{eq:normalization-Rhat-bounds}
\end{align}
Each estimate is understood over the corresponding index range.
The constant $C$ is independent of $m,j,K,H$.
\end{proposition}

\begin{proof}
We use the previously established quantitative graded expansion
lemma. All decompositions into principal parts and remainders
are made using its fixed frequency cutoffs.

\medskip
\noindent
\textbf{Step 1. The graded expansion and triangular dependence.}

The graded expansion reads
\begin{align}
&(I+\mathcal Q_m)^*(P_0+V)(I+\mathcal Q_m)
\nonumber\\
&=
P_0+V+\sum_{r=1}^m[Q'_r,\Delta]
+\sum_{r=1}^m(Q''_rP_0+P_0Q''_r)
\nonumber\\
&\quad+
\frac12\sum_{j=1}^{2m+1}(S_jP_0+P_0S_j)
+\frac12\sum_{j=1}^{2m+1}(R_jP_0+P_0R_j)
\nonumber\\
&\quad+
\sum_{j=2}^{m+1}\widetilde S_j
+\sum_{j=2}^{2m+3}(T_j+U_j).
\label{eq:normalization-graded}
\end{align}
For each $r\ge2$, the coefficients $\widetilde S_r,T_r$
depend only on $Q'_i,Q''_i$ with $i<r$.
The coefficient $S_r$ depends on $Q'_i$ with $i\le r$
and on $Q''_i$ with $i<r$.

Consequently, the recursive construction has the order
\[
Q'_1\longrightarrow S_1\longrightarrow Q''_1
\longrightarrow Q'_2\longrightarrow S_2\longrightarrow Q''_2
\longrightarrow\cdots.
\]
Because the decomposition rules are fixed, the already constructed
coefficients do not change when the normalization order is increased.

\medskip
\noindent
\textbf{Step 2. Small-gap off-block terms.}

For a self-adjoint operator
\[
F\in\widetilde{\mathcal L}_\tau^{-(r-1)}(1),
\]
decompose
\[
F=F_{\mathrm D}
+F_{\mathrm{ND,large}}
+F_{\mathrm{ND,small}},
\]
where the large-gap part is restricted to pairs in distinct
Bourgain blocks satisfying
\[
\bigl||n|^2-|n'|^2\bigr|
>
\frac14(|n|+|n'|)^{\tau_0}.
\]
All three restrictions preserve self-adjointness.

On the small-gap off-block support, the block separation property
implies, outside a bounded low-frequency region,
\[
|n-n'|
\ge
c\bigl(1+\max(|n|,|n'|)\bigr)^{\tau_0}.
\]
Thus
\begin{align*}
e^{-\rho\phi(|n-n'|)}
&\le
C\exp\!\left[
-c_\rho
\bigl(\log(1+\max(|n|,|n'|))\bigr)^a
\right]\\
&\le
C\exp\!\left(c_1N^{\frac a{a-1}}\right)
\bigl(1+\max(|n|,|n'|)\bigr)^{-\tau N}
\end{align*}
for every integer $N\ge0$.
On the other hand,
\[
\frac{W_{q+r-1+N}}{W_{q+r-1}}
\ge
\exp\!\left(\kappa N^{\frac a{a-1}}\right).
\]
Choosing $\kappa\ge c_1$ in the initial parameter selection gives
\[
F_{\mathrm{ND,small}}\in\mathcal R_\tau^{-(r-1)},
\qquad
|F_{\mathrm{ND,small}}|_{r-1}
\le C_{\mathrm g}\|F\|_{r-1,1}.
\]
The bounded low-frequency region is absorbed into the fixed
constants. In particular, the constants are independent of $r$.

\medskip
\noindent
\textbf{Step 3. Construction of $Q'_1$.}

On the large-gap off-block support, define
\[
\Pi_nQ'_1\Pi_{n'}
=
\frac{\Pi_n(Q_V)_{\mathrm{ND,large}}\Pi_{n'}}
{|n'|^2-|n|^2},
\]
and set all other matrix entries equal to zero.
Since
\[
\Delta e^{in\cdot x}=-|n|^2e^{in\cdot x},
\]
we have
\[
[Q'_1,\Delta]=-(Q_V)_{\mathrm{ND,large}}.
\]
The numerator is self-adjoint, whereas the denominator changes
sign under interchange of $n,n'$. Hence
\[
(Q'_1)^*=-Q'_1.
\]

At grade zero the two temporal indices coincide:
\[
(0+0-1)_+=(0+1-1)_+=0.
\]
Therefore
\[
\|Q_V\|_{0,1}=\|Q_V\|_{0,0}\le H.
\]
The homological equation estimate yields
\[
\|Q'_1\|_{1,0}\le C_{\mathrm h}H,
\qquad
\|[\Delta,Q'_1]\|_{0,1}\le H.
\]
Thus \eqref{eq:normalization-prime-bounds} holds at $r=1$
provided
\[
\sqrt K\ge C_{\mathrm h}.
\]

\medskip
\noindent
\textbf{Step 4. Explicit construction of $S_r$ and $Q''_r$.}

Suppose that $Q'_r$ has been constructed and that
$Q''_1,\ldots,Q''_{r-1}$ are already known.
Define
\begin{align}
\Pi_nS_r\Pi_{n'}
&=
\mathbf1_{\left\{
|n-n'|\le
\frac{\max(|n|,|n'|)}{10(r+3)}
\right\}}
\nonumber\\
&\quad\times
\Pi_n\Bigg[
\sum_{\substack{i+j=r+1\\i,j\ge1}}(Q'_i)^*Q'_j+
\sum_{\substack{i+j=r\\i,j\ge1}}
(Q'_i)^*Q''_j+(Q''_j)^*Q'_i
\nonumber\\
&\qquad\quad+
\sum_{\substack{i+j=r-1\\i,j\ge1}}
(Q''_i)^*Q''_j
\Bigg]\Pi_{n'}.
\label{eq:normalization-S-construction}
\end{align}

In the first sum, the indices do not exceed $r$.
In the second sum, the double-prime index does not exceed $r-1$.
In the third sum, both indices do not exceed $r-2$.
Thus every operator on the right-hand side is already known;
in particular, $Q''_r$ does not occur there.

Set
\[
Q''_r=-\frac12S_r.
\]
At the first step this reduces to
\[
\Pi_nS_1\Pi_{n'}
=
\mathbf1_{\left\{
|n-n'|\le\frac{\max(|n|,|n'|)}{40}
\right\}}
\Pi_n(Q'_1)^*Q'_1\Pi_{n'},
\]
so the recursion starts without any double-prime input.

Symmetric summation and the symmetric cutoff imply
$S_r^*=S_r$. The ordinary composition estimate and
\[
\sum_{\ell=1}^{r}
\frac1{\ell^2(r+1-\ell)^2}
\le\frac{C}{(r+1)^2}
\]
give
\[
\|S_r\|_{r+1,0}
\le C_1\frac{K^rH^{r+1}}{(r+1)^2}.
\]

To estimate the commutator, apply the product rule to each term.
For example,
\[
[\Delta,(Q'_i)^*Q'_j]
=
[\Delta,(Q'_i)^*]Q'_j
+(Q'_i)^*[\Delta,Q'_j].
\]
Every factor with flag $(1)$ is paired with a positive-order
factor with flag $(0)$.
The mixed composition theorem therefore applies, including the
endpoint where the commutator has order zero.

The total smoothing order after taking the commutator is $r\tau$.
The corresponding composition cutoff has denominator $10(r+2)$,
whereas the cutoff in
\eqref{eq:normalization-S-construction} has denominator $10(r+3)$.
Since the latter support is narrower, it eliminates the
complementary output-cutoff remainders.
Also, taking the commutator with $\Delta$ commutes with these
matrix cutoffs. Consequently,
\[
[\Delta,S_r]\in\widetilde{\mathcal L}_\tau^{-r}(0),
\]
and
\[
\|[\Delta,S_r]\|_{r,0}
\le C_1\frac{K^{r-\frac12}H^{r+1}}{(r+1)^2}.
\]

It follows that
\[
\|Q''_r\|_{r+1,0}
\le\frac{C_1}{2}\frac{K^rH^{r+1}}{(r+1)^2},
\]
and
\[
\|[\Delta,Q''_r]\|_{r,0}
\le\frac{C_1}{2}
\frac{K^{r-\frac12}H^{r+1}}{(r+1)^2}.
\]
These imply \eqref{eq:normalization-double-bounds} whenever
\[
\sqrt K\ge C_1/2.
\]
Moreover, the definition gives the exact cancellation
\[
Q''_rP_0+P_0Q''_r
+\frac12(S_rP_0+P_0S_r)=0.
\]

\medskip
\noindent
\textbf{Step 5. Construction of $Q'_r$ for $r\ge2$.}

Assume that all operators of index less than $r$ have been
constructed and satisfy
\eqref{eq:normalization-prime-bounds}--%
\eqref{eq:normalization-double-bounds}.
By triangular dependence, $\widetilde S_r,T_r$ are already
determined. The quantitative graded expansion gives
\[
\|\widetilde S_r\|_{r-1,1}
\le C_0\frac{K^{r-\frac32}H^{r-1}}{r^2},
\]
and
\[
\|T_r\|_{r-1,0}
\le C_0\frac{K^{r-\frac32}H^r}{r^2}.
\]
Since $H\ge1$ and the flag-$(0)$ class embeds into the
flag-$(1)$ class, increasing $C_0$ if necessary yields
\[
\|\widetilde S_r+T_r\|_{r-1,1}
\le C_0\frac{K^{r-\frac32}H^r}{r^2}.
\]

Define
\[
\Pi_nQ'_r\Pi_{n'}
=
\frac{
\Pi_n(\widetilde S_r+T_r)_{\mathrm{ND,large}}\Pi_{n'}
}{
|n'|^2-|n|^2
}
\]
on the large-gap off-block support, and set the remaining
entries equal to zero. Then
\[
(Q'_r)^*=-Q'_r,
\qquad
[Q'_r,\Delta]
=-(\widetilde S_r+T_r)_{\mathrm{ND,large}}.
\]

The source class
$\widetilde{\mathcal L}_\tau^{-(r-1)}(1)$ has temporal weight
$W_{q+r-1}$, which is precisely the temporal weight of
$\mathcal L_\tau^{-r}(0)$.
Division by the large spectral gap supplies at least $\tau$
additional spatial smoothing, since $\tau\le\tau_0$.
The source cutoff also agrees with the support requirement for
$\mathcal L_\tau^{-r}(0)$.
Hence
\[
\|Q'_r\|_{r,0}
\le C_{\mathrm h}C_0
\frac{K^{r-\frac32}H^r}{r^2}.
\]
The commutator identity directly gives
\[
\|[\Delta,Q'_r]\|_{r-1,1}
\le C_0\frac{K^{r-\frac32}H^r}{r^2}.
\]
Thus \eqref{eq:normalization-prime-bounds} follows if
\[
K\ge C_{\mathrm h}C_0,
\qquad
\sqrt K\ge C_0.
\]

We then construct $S_r,Q''_r$ as in Step 4.
One fixed choice of $K$ satisfying all the preceding inequalities
closes the induction at every order.

\medskip
\noindent
\textbf{Step 6. The block-diagonal potential and homological
remainders.}

Define
\[
V^{(m)}
=
(Q_V)_{\mathrm D}
+\sum_{r=2}^{m}(\widetilde S_r+T_r)_{\mathrm D},
\]
and
\[
\widehat R_1
=
R_V+(Q_V)_{\mathrm{ND,small}},
\]
\[
\widehat R_r
=
(\widetilde S_r+T_r)_{\mathrm{ND,small}},
\qquad 2\le r\le m.
\]
These operators are self-adjoint.
The operator $V^{(m)}$ is block diagonal, and therefore
\[
[V^{(m)},\widetilde\Delta]=0.
\]
Step 2 gives
\[
|\widehat R_1|_0\le CH,
\]
and
\[
|\widehat R_r|_{r-1}
\le C\frac{K^{r-\frac32}H^r}{r^2},
\qquad r\ge2.
\]

The homological equations imply
\[
V+[Q'_1,\Delta]
=(Q_V)_{\mathrm D}+\widehat R_1,
\]
and
\[
\widetilde S_r+T_r+[Q'_r,\Delta]
=
(\widetilde S_r+T_r)_{\mathrm D}+\widehat R_r,
\qquad r\ge2.
\]
Summing these identities gives
\begin{align*}
V+\sum_{r=1}^{m}[Q'_r,\Delta]
+\sum_{r=2}^{m}(\widetilde S_r+T_r)
=
V^{(m)}+\sum_{r=1}^{m}\widehat R_r.
\end{align*}

The terms $\widehat R_r$ are distinct from the terms $U_r$:
the former arise from small-gap off-block restrictions, whereas
the latter arise from cutoff and composition remainders,
including products containing $R_V$.

\medskip
\noindent
\textbf{Step 7. Completion of the normal form.}

Substitute the preceding identity and
\[
Q''_rP_0+P_0Q''_r
+\frac12(S_rP_0+P_0S_r)=0,
\qquad 1\le r\le m,
\]
into \eqref{eq:normalization-graded}.
This gives exactly \eqref{eq:normalization-final}.

The remaining class memberships and estimates
\eqref{eq:normalization-SR-bounds}--%
\eqref{eq:normalization-TU-bounds}
follow from the quantitative graded expansion lemma.
The estimates for $\widehat R_j$ were proved in Step 6.

Since $K,H$ and all calculus parameters are independent of $m$,
these estimates imply exponential bounds in the grading indices
with a base independent of $m$.
Finally, for each fixed $m$, $V^{(m)}$ is a finite sum of operators
uniformly bounded in time on $L^2$, so it has the asserted
boundedness.
\end{proof}

\begin{lemma}[ Finite-order approximate inversion] \label{L4.1}
Given $Q_j\in\mathcal{L}_\tau^{-j}(0)$, $1\le j\le m$, there exist
\[
 P_j\in\mathcal{L}_\tau^{-j}(0)\ (1\le j\le m),\quad
 Z_j\in\mathcal{L}_\tau^{-j}(0)\ (m+1\le j\le2m),\quad
 E_j\in\mathcal{R}_j^{-\infty}\ (2\le j\le2m),
\]
such that
\begin{equation}\label{eq:inverse}
 (I+\sum_{j=1}^mQ_j)(I+\sum_{j=1}^mP_j)
 =I+\sum_{j=m+1}^{2m}Z_j+\sum_{j=2}^{2m}E_j.
\end{equation}
Set $a_j=\|Q_j\|_{j,0}$. For a suitable fixed $C_2\ge1$,
\begin{align}
 \|P_j\|_{j,0}
 &\le\sum_{\ell=1}^j\ \sum_{\substack{i_1+\cdots+i_\ell=j\\1\le i_\nu\le m}}
 C_2^{\ell-1}\prod_{\nu=1}^\ell a_{i_\nu},\label{eq:inverseP}\\
 \|Z_j\|_{j,0}
 &\le\sum_{\ell=2}^j\ \sum_{\substack{i_1+\cdots+i_\ell=j\\1\le i_\nu\le m}}
 C_2^{\ell-1}\prod_{\nu=1}^\ell a_{i_\nu},\label{eq:inverseZ}\\
 |E_j|_j
 &\le\sum_{\ell=2}^j\ \sum_{\substack{i_1+\cdots+i_\ell=j\\1\le i_\nu\le m}}
 C_2^{\ell-1}\prod_{\nu=1}^\ell a_{i_\nu}.\label{eq:inverseE}
\end{align}
Each inequality holds over its corresponding index range.\ In particular, if $a_j\le C_0^{j+1}$, all these seminorms are bounded by $C_3^{j+1}$, with $C_3$ independent of $m$.
\end{lemma}
\begin{proof}
Take $P_1=-Q_1$. For $2\le j\le m$, define recursively
\[
 P_j=-Q_j-\sum_{i+k=j}\mathcal{M}(Q_i,P_k).
\]
The remaining sum $\sum_{i+k=j}\mathcal{E}(Q_i,P_k)$ defines $E_j$.
For $j>m$, split $\sum_{i+k=j}Q_iP_k$ into $Z_j+E_j$ by Proposition 2.7.
Expanding the product gives~\eqref{eq:inverse}. The seminorm recursion satisfies
$p_j\le a_j+C_2\sum_{i+k=j}a_ip_k$. Expansion over ordered integer compositions gives~\eqref{eq:inverseP}; substitution into the principal and remainder estimates gives the other two inequalities.
\end{proof}

\begin{lemma}[Energy estimates for normalization remainders]\label{L4.2}Let
\[
p=\frac{\sigma}{\sigma-1},
\qquad
P_0=i\partial_t-\Delta,
\]
and assume that the fixed parameters satisfy $0<\tau<1$.\ Let $u$ solve
\[
(P_0+V)u=0,
\qquad
u(0)=u_0.
\]
Given $N\ge1$, choose $m$ such that
\begin{equation}\label{eq:mN}
 N+3\le(m+2)\tau<N+4.
\end{equation}
Regroup the corrections by smoothing order as
\[
 Q_1=Q'_1,\qquad Q_j=Q'_j+Q''_{j-1}\ (2\le j\le m),\qquad
 Q_{m+1}=Q''_m.
\]
Apply Lemma \ref{L4.1} to these $m+1$ terms, and write
\[
 T_m=I+\mathcal Q_m,\qquad B_m=I+\sum_{j=1}^{m+1}P_j,
 \qquad T_mB_m=I+E^{\mathrm{inv}}_m.
\]
Let $\mathcal D_m$ denote all terms in (\ref{2.14*}) other than $P_0+V^{(m)}$.
Set $v=B_mu$. Then
\begin{equation}\label{eq:fg}
 (P_0+V^{(m)})v=f+g,\qquad
 f=-\mathcal D_mv,\qquad
 g=T_m^*[P_0+V,E^{\mathrm{inv}}_m]u.
\end{equation}
There exists $C$ independent of $m,N,t$ such that
\begin{equation}\label{eq:fgbound}
 \|f(t)\|_{\widetilde H^N}+ \|g(t)\|_{\widetilde H^N}
 \le e^{C(m+1)^p}\|u_0\|_2.
\end{equation}
\end{lemma}

\begin{proof}
We first work with smooth initial data. Since $v=B_mu$,
\[
T_mv=(I+E_m^{\mathrm{inv}})u.
\]
Using $(P_0+V)u=0$, we obtain
\begin{align*}
T_m^*(P_0+V)T_mv
&=
T_m^*(P_0+V)(I+E_m^{\mathrm{inv}})u\\
&=
T_m^*[P_0+V,E_m^{\mathrm{inv}}]u.
\end{align*}
The normalization identity therefore gives
\[
(P_0+V^{(m)})v+\mathcal D_mv
=
T_m^*[P_0+V,E_m^{\mathrm{inv}}]u.
\]
This proves the asserted equation and the formulas for $f,g$.\\
The corrected normalization proposition and Lemma \ref{L4.1} give
exponential bounds in the grading indices for all coefficient
seminorms. By construction,
\[
B_m=I+\sum_{j=1}^{m+1}P_j,
\qquad
P_j\in\mathcal L_\tau^{-j}(0).
\]
Lemma \ref{L4.1} gives a constant $C_3\ge1$, independent of $m,j$,
such that
\[
\|P_j\|_{j,0}\le C_3^{j+1}.
\]
Since $j\ge1$, the defining estimate for this class is
\begin{align*}
\sup_t\|\Pi_n\partial_t^qP_j(t)\Pi_{n'}\|
&\le
C_3^{j+1}W_{q+j-1}
e^{-\rho\phi(|n-n'|)}
\langle n-n'\rangle^{-(d+2)}\\
&\quad\times
(1+\max(|n|,|n'|))^{-j\tau}
\mathbf1_{\left\{
|n-n'|\le
\frac{\max(|n|,|n'|)}{10(1+j)}
\right\}}.
\end{align*}
On this support,
\[
\frac9{10}
\le\frac{\langle n\rangle}{\langle n'\rangle}
\le\frac{10}{9},
\qquad
\frac{\langle n\rangle^{j\tau}}
{(1+\max(|n|,|n'|))^{j\tau}}
\le1.
\]
In view of proposition \ref{2.5}, one has
\[
\sup_t
\|\partial_t^qP_j(t)\|_{H^s\to H^{s+j\tau}}
\le
C_d\left(\frac{10}{9}\right)^{|s|}
C_3^{j+1}W_{q+j-1}.
\]
Taking $s=0$, $q=0$, and using
$H^{j\tau}\hookrightarrow L^2$, we obtain
\[
\sup_t\|P_j(t)\|_{L^2\to L^2}
\le C_dC_3^{j+1}W_{j-1}.
\]
Taking $s=-2$, $q=0$, and using
$H^{-2+j\tau}\hookrightarrow H^{-2}$ gives
\[
\sup_t\|P_j(t)\|_{H^{-2}\to H^{-2}}
\le
C_d\left(\frac{10}{9}\right)^2C_3^{j+1}W_{j-1}.
\]
Finally, taking $s=0$, $q=1$, and using
$H^{j\tau}\hookrightarrow H^{-2}$ gives
\[
\sup_t\|\partial_tP_j(t)\|_{L^2\to H^{-2}}
\le C_dC_3^{j+1}W_j.
\]
Consequently, the triangle inequality implies
\begin{align*}
&\sup_t\left(
\|B_m(t)\|_{L^2\to L^2}
+\|B_m(t)\|_{H^{-2}\to H^{-2}}
+\|\partial_tB_m(t)\|_{L^2\to H^{-2}}
\right)\\
&\qquad\le
2+C\sum_{j=1}^{m+1}C_3^{j+1}(W_{j-1}+W_j).
\end{align*}
Here the constant $2$ accounts for the identity operator in the
first two norms.\\
Since
\[
W_j=M^j(j!)^2e^{\kappa j^p},
\qquad p>1,
\]
we have
\begin{align*}
\log(C_3^{j+1}W_j)
&=(j+1)\log C_3+j\log M+2\log(j!)+\kappa j^p \le C_4(j+1)^p,
\end{align*}
where we used
\[
j\log(j+1)\le C_p(j+1)^p.
\]
Thus
\begin{align}\nonumber&\|B_m\|_{L^2\to L^2}
+\|B_m\|_{H^{-2}\to H^{-2}}
+\|\partial_tB_m\|_{L^2\to H^{-2}}\\
\nonumber&\leq 2+C\sum_{j=1}^{m+1}C_3^{j+1}(W_{j-1}+W_j)\\
\nonumber&\le
2+C(m+1)e^{C_4(m+2)^p}\\
&\label{30m}\le e^{C(m+1)^p}.
\end{align}
Since $V$ is real-valued, the evolution is unitary on $L^2$:
\[
\|u(t)\|_{L^2}=\|u_0\|_{L^2}.
\]
The equation
\[
i\partial_tu=\Delta u-Vu
\]
and the uniform $L^2$ boundedness of multiplication by $V$ imply
\[
\|\partial_tu(t)\|_{H^{-2}}
\le C\|u_0\|_{L^2}.
\]
By (\ref{30m}), we have
\[
\|v(t)\|_{L^2}
\le e^{C(m+1)^p}\|u_0\|_{L^2}.
\]
Using
\[
\partial_tv=(\partial_tB_m)u+B_m\partial_tu,
\]
we also obtain
\[
\|\partial_tv(t)\|_{H^{-2}}
\le e^{C(m+1)^p}\|u_0\|_{L^2}.
\]
Since $\Delta:L^2\to H^{-2}$ is bounded,
\[
\|P_{0}v(t)\|_{H^{-2}}+\|v(t)\|_{L^2}
\le e^{C(m+1)^p}\|u_0\|_{L^2}.
\]
The normalization remainder is
\begin{align*}
\mathcal D_m
&=
\frac12\sum_{j=m+1}^{2m+1}(S_jP_0+P_0S_j)\\
&\quad+
\frac12\sum_{j=1}^{2m+1}(R_jP_0+P_0R_j)
+\widetilde S_{m+1}\\
&\quad+
\sum_{j=m+1}^{2m+3}T_j
+\sum_{j=2}^{2m+3}U_j
+\sum_{j=1}^{m}\widehat R_j.
\end{align*}
For $j\ge m+1$,
\[
S_j\in\widetilde{\mathcal L}_\tau^{-(j+1)}(0),
\qquad
(j+1)\tau\ge(m+2)\tau\ge N+3.
\]
Thus $S_j:H^{-2}\to H^N$ is bounded, with the bound from Step~2.
It follows that
\[
\|S_jP_0v\|_{H^N}
\le e^{C(m+1)^p}\|u_0\|_{L^2}.
\]
For the other ordering, write
\[
P_0S_jv
=
S_jP_0v+i(\partial_tS_j)v-[\Delta,S_j]v.
\]
Time differentiation preserves the spatial smoothing order, so
\[
\|(\partial_tS_j)v\|_{H^N}
\le e^{C(m+1)^p}\|u_0\|_{L^2}.
\]
The corrected normalization proposition gives
\[
[\Delta,S_j]\in
\widetilde{\mathcal L}_\tau^{-j}(0).
\]
Since
\[
j\tau\ge(m+1)\tau\ge N+3-\tau>N,
\]
we obtain
\[
\|[\Delta,S_j]v\|_{H^N}
\le e^{C(m+1)^p}\|u_0\|_{L^2}.
\]
This controls all terms involving $S_j$.\ For $R_j$, infinite smoothing gives
\[
R_j:H^{-2}\to H^N,
\qquad
R_j:L^2\to H^{N+2},
\qquad
\partial_tR_j:L^2\to H^N,
\]
with mapping norms bounded by $e^{C(m+1)^p}$. Using
\[
P_0R_jv
=
i(\partial_tR_j)v+iR_j\partial_tv-\Delta R_jv,
\]
we conclude that
\[
\|R_jP_0v\|_{H^N}
+\|P_0R_jv\|_{H^N}
\le e^{C(m+1)^p}\|u_0\|_{L^2}.
\]

Next,
\[
\widetilde S_{m+1}\in
\widetilde{\mathcal L}_\tau^{-m}(1),
\qquad
T_j\in
\widetilde{\mathcal L}_\tau^{-(j-1)}(0),
\quad j\ge m+1.
\]
Their smoothing orders are at least $m\tau$, and
\[
m\tau\ge N+3-2\tau>N.
\]
Therefore
\[
\|\widetilde S_{m+1}v\|_{H^N}
+
\sum_{j=m+1}^{2m+3}\|T_jv\|_{H^N}
\le e^{C(m+1)^p}\|u_0\|_{L^2}.
\]
The flag $(1)$ of $\widetilde S_{m+1}$ changes its weight index,
but does not reduce its spatial smoothing order.

Finally, the infinite smoothing of $U_j,\widehat R_j$ gives
\[
\sum_{j=2}^{2m+3}\|U_jv\|_{H^N}
+
\sum_{j=1}^{m}\|\widehat R_jv\|_{H^N}
\le e^{C(m+1)^p}\|u_0\|_{L^2}.
\]
Combining these estimates yields
\[
\|f(t)\|_{H^N}
\le e^{C(m+1)^p}\|u_0\|_{L^2}.
\]
Applying Lemma \ref{L4.1} to $m+1$ terms gives
\[
E_m^{\mathrm{inv}}
=
\sum_{j=m+2}^{2m+2}Z_j
+\sum_{j=2}^{2m+2}E_j,
\]
where
\[
Z_j\in\mathcal L_\tau^{-j}(0),
\qquad
E_j\in\mathcal R_j^{-\infty}.
\]
The finite-order part has smoothing order at least
\[
(m+2)\tau\ge N+3.
\]
Consequently,
\begin{align*}
&\|E_m^{\mathrm{inv}}\|_{H^{-2}\to H^N}
+\|E_m^{\mathrm{inv}}\|_{L^2\to H^{N+2}}+
\|\partial_tE_m^{\mathrm{inv}}\|_{L^2\to H^N}
\le e^{C(m+1)^p}.
\end{align*}
For the Laplacian commutator, use
\[
[\Delta,E_m^{\mathrm{inv}}]u
=
\Delta E_m^{\mathrm{inv}}u
-E_m^{\mathrm{inv}}\Delta u.
\]
Since
\[
\|\Delta u\|_{H^{-2}}\le C\|u\|_{L^2},
\]
we have
\begin{align*}
\|[\Delta,E_m^{\mathrm{inv}}]u\|_{H^N}
&\le
C\|E_m^{\mathrm{inv}}u\|_{H^{N+2}}
+\|E_m^{\mathrm{inv}}\Delta u\|_{H^N}\le e^{C(m+1)^p}\|u_0\|_{L^2}.
\end{align*}

We below estimate multiplication by $V$ on $H^N$.
The Fourier convolution formula and
\[
\langle n\rangle^N
\le C^N\langle n-k\rangle^N\langle k\rangle^N
\]
give
\[
\|V(t)w\|_{H^N}
\le
C^N
\left(\sum_{k\in\mathbb Z^d}
\langle k\rangle^N|v_k(t)|\right)
\|w\|_{H^N}.
\]
The assumed spatial decay implies
\[
\sup_t\sum_k\langle k\rangle^N|v_k(t)|
\le e^{C(N+1)^p}.
\]
Indeed, reserve half of the spatial decay for summability and
use the other half to absorb the polynomial weight, applying
\[
\sup_{x\ge0}e^{Nx-cx^a}
\le e^{CN^p}.
\]
Thus
\[
\sup_t\|V(t)\|_{H^N\to H^N}
\le e^{C(N+1)^p}.
\]
Since $N=O(m+1)$, Step~5 and the uniform $L^2$ boundedness of $V$
yield
\begin{align*}
\|[V,E_m^{\mathrm{inv}}]u\|_{H^N}
&\le
\|VE_m^{\mathrm{inv}}u\|_{H^N}
+\|E_m^{\mathrm{inv}}Vu\|_{H^N}\le e^{C(m+1)^p}\|u_0\|_{L^2}.
\end{align*}
The commutator identity
\[
[P_0+V,E_m^{\mathrm{inv}}]
=
i\partial_tE_m^{\mathrm{inv}}
-[\Delta,E_m^{\mathrm{inv}}]
+[V,E_m^{\mathrm{inv}}]
\]
and Steps~5--6 imply
\[
\|[P_0+V,E_m^{\mathrm{inv}}]u\|_{H^N}
\le e^{C(m+1)^p}\|u_0\|_{L^2}.
\]
Moreover, the mapping estimates for the normalizing corrections
give
\[
\|T_m^*\|_{H^N\to H^N}
\le e^{C(m+1)^p}.
\]
Hence
\[
\|g(t)\|_{H^N}
\le e^{C(m+1)^p}\|u_0\|_{L^2}.
\]
The frequency comparability in the block construction gives
\[
\|w\|_{\widetilde H^N}
\le C_{\mathrm{eq}}^{\,N}\|w\|_{H^N}.
\]
Since $N=O(m+1)$, this factor is absorbed into
$e^{C(m+1)^p}$. Combining the estimates for $f$ and $g$ proves
\[
\|f(t)\|_{\widetilde H^N}
+\|g(t)\|_{\widetilde H^N}
\le e^{C(m+1)^p}\|u_0\|_{L^2}.
\]

All bounds depend on the initial data only through its $L^2$ norm.
Applying the estimates to differences of smooth approximating
solutions therefore extends the conclusion to $L^2$ initial data,
with the transformed equation interpreted in the distributional
sense.
\end{proof}

\section{Proof of Theorem \ref{thm:main}}
\begin{proof}
We use the corrected normalization proposition, the
approximate-inverse lemma, and Lemma \ref{L4.2}, with all parameters
fixed independently of the normalization order.

The spatial Fourier decay implies that, for every fixed $s\ge0$,
\[
\sup_t\sum_{k\in\mathbb Z^d}
\langle k\rangle^s|v_k(t)|<\infty.
\]
The Fourier convolution estimate therefore gives
\[
\sup_t\|V(t)\|_{H^s\to H^s}<\infty.
\]
The time regularity assumption also ensures continuity of
$t\mapsto V(t)$ in the required operator topology. Denote the propagator by $U(t,0): u(t)=U(t,0)u_0.$ Since $V$ is real-valued, the propagator is unitary on $L^2$:
\[
\|U(t,0)\|_{L^2\to L^2}=1.
\]
Let $N\ge1$ be an integer. Choose $m$ such that
\[
N+3\le(m+2)\tau<N+4,
\]
where $0<\tau<1$ is fixed. In particular,
\[
m+1\le C_\tau(N+1).
\]
As in Lemma \ref{L4.2}, let
\[
T_m=I+\mathcal Q_m,
\qquad
B_m=I+\sum_{j=1}^{m+1}P_j,
\qquad
T_mB_m=I+E_m^{\mathrm{inv}}.
\]
For
\[
v=B_mu,
\]
the transformed equation is
\[
(i\partial_t-\Delta+V^{(m)})v=f+g,
\]
where
\[
\|f(t)\|_{\widetilde H^N}
+\|g(t)\|_{\widetilde H^N}
\le
e^{C(m+1)^p}\|u_0\|_{L^2}.
\]
Define
\[
\widetilde\Lambda^N=(1-\widetilde\Delta)^{N/2}.
\]
Since $\Delta$ and $\widetilde\Delta$ are diagonal in the Fourier
basis,
\[
[\Delta,\widetilde\Lambda^N]=0.
\]
Furthermore, $V^{(m)}$ is block diagonal and $\widetilde\Delta$
is scalar on each block. Hence
\[
[V^{(m)},\widetilde\Lambda^N]=0.
\]
Applying $\widetilde\Lambda^N$ to the transformed equation yields
\[
i\partial_t(\widetilde\Lambda^Nv)
=
(\Delta-V^{(m)})\widetilde\Lambda^Nv
+\widetilde\Lambda^N(f+g).
\]
The operator $\Delta-V^{(m)}(t)$ is self-adjoint.
Consequently, the homogeneous part contributes nothing to the
derivative of the squared $L^2$ norm, and
\[
\left|
\frac{d}{dt}\|\widetilde\Lambda^Nv\|_{L^2}^2
\right|
\le
2\|\widetilde\Lambda^Nv\|_{L^2}
\|\widetilde\Lambda^N(f+g)\|_{L^2}.
\]
The standard inhomogeneous energy inequality gives
\begin{align*}
\|v(t)\|_{\widetilde H^N}
&\le
\|v(0)\|_{\widetilde H^N}+
\int_{\min(0,t)}^{\max(0,t)}
\|f(\sigma)+g(\sigma)\|_{\widetilde H^N}\,d\sigma.
\end{align*}
Thus
\[
\|v(t)\|_{\widetilde H^N}
\le
\|v(0)\|_{\widetilde H^N}
+
|t|e^{C(m+1)^p}\|u_0\|_{L^2}.
\]
The approximate-inverse lemma gives
\[
P_j\in\mathcal L_\tau^{-j}(0),
\qquad
\|P_j\|_{j,0}\le C_3^{j+1},
\]
with $C_3$ independent of $m,j$.
The mapping estimate implies
\[
\|P_j\|_{H^N\to H^N}
\le
C_d\left(\frac{10}{9}\right)^N
C_3^{j+1}W_{j-1}.
\]
Using
\[
\log W_b\le C(b+1)^p,
\qquad N=O(m+1),
\]
and summing over $1\le j\le m+1$, we obtain
\[
\sup_t\|B_m(t)\|_{H^N\to H^N}
\le e^{C(m+1)^p}.
\]
The ordinary and block Sobolev norms satisfy
\[
C_{\mathrm{eq}}^{-N}\|w\|_{H^N}
\le
\|w\|_{\widetilde H^N}
\le
C_{\mathrm{eq}}^N\|w\|_{H^N},
\]
where $C_{\mathrm{eq}}$ is independent of $N$.
It follows that
\[
\|v(0)\|_{\widetilde H^N}
=
\|B_m(0)u_0\|_{\widetilde H^N}
\le
e^{C(m+1)^p}\|u_0\|_{H^N}.
\]
Therefore
\[
\|v(t)\|_{\widetilde H^N}
\le
e^{C(m+1)^p}(1+|t|)\|u_0\|_{H^N}.
\]
By
\[
T_mv=T_mB_mu=u+E_m^{\mathrm{inv}}u,
\]
we have
\[
u=T_mv-E_m^{\mathrm{inv}}u.
\]
The normalizing corrections satisfy
\[
\|T_m\|_{H^N\to H^N}
\le e^{C(m+1)^p}.
\]
The finite-order part of $E_m^{\mathrm{inv}}$ has smoothing order
at least
\[
(m+2)\tau\ge N+3,
\]
and its remaining part is infinitely smoothing. Hence
\[
\|E_m^{\mathrm{inv}}\|_{L^2\to H^N}
\le e^{C(m+1)^p}.
\]
Using $L^2$ conservation and the preceding bound for $v$, we obtain
\begin{align*}
\|u(t)\|_{H^N}
&\le
e^{C(m+1)^p}\|v(t)\|_{H^N}
+
e^{C(m+1)^p}\|u_0\|_{L^2}\le
e^{C(m+1)^p}(2+|t|)\|u_0\|_{H^N}.
\end{align*}
Since $m+1\le C_\tau(N+1)$ and $N\ge1$, this yields
\[
\|U(t,0)\|_{H^N\to H^N}
\le e^{CN^p}(2+|t|),
\qquad N\in\mathbb N.
\]

Fix $s>0$ and choose an integer $N\ge s$.
Interpolating the operator bounds
\[
\|U(t,0)\|_{L^2\to L^2}=1
\]
and
\[
\|U(t,0)\|_{H^N\to H^N}
\le e^{CN^p}(2+|t|)
\]
gives
\begin{align*}
\|U(t,0)\|_{H^s\to H^s}
&\le
\|U(t,0)\|_{H^N\to H^N}^{s/N}\le
\exp\!\left(
CsN^{p-1}
+\frac{s\log(2+|t|)}{N}
\right).
\end{align*}
Choose
\[
N=
\max\left\{
1,s,[\log(2+|t|)]^{1/p}
\right\}.
\]
For sufficiently large $|t|$,
\[
N\asymp[\log(2+|t|)]^{1/p}.
\]
Consequently,
\[
N^{p-1}
\le C[\log(2+|t|)]^{(p-1)/p},
\]
and
\[
\frac{\log(2+|t|)}{N}
\le C[\log(2+|t|)]^{(p-1)/p}.
\]
Since
\[
\frac{p-1}{p}=\frac1\sigma,
\]
the interpolation estimate becomes
\[
\|U(t,0)\|_{H^s\to H^s}
\le
\exp\!\left(
c_s[\log(2+|t|)]^{1/\sigma}
\right)
\]
for sufficiently large $|t|$.
Absorbing bounded time intervals into a multiplicative constant
$C_s$ proves
\[
\|u(t)\|_{H^s}
\le
C_s\exp\!\left(
c_s[\log(2+|t|)]^{1/\sigma}
\right)\|u_0\|_{H^s},
\qquad t\in\mathbb R.
\]
\end{proof}

\section{Appendix}
\begin{lemma}\label{t1}[Temporal convolution estimate]
For nonnegative integers $q,b_1,b_2$, we have
\begin{equation}\label{eq:p27-time}
\sum_{\ell=0}^q
\binom q\ell W_{\ell+b_1}W_{q-\ell+b_2}
\le 3W_{q+b_1+b_2}.
\end{equation}
\end{lemma}
\begin{proof}
Put
\[
L=q+b_1+b_2,\qquad r=\ell+b_1.
\]
Since $p>1$,
\[
r^p+(L-r)^p\le L^p.
\]
Consequently,
\[
\frac{\binom q\ell W_rW_{L-r}}{W_L}
\le
\frac{\binom q\ell}{\binom Lr^2}
\le
\binom Lr^{-1},
\]
where
\[
\binom q\ell\le\binom L{\ell+b_1}.
\]
Summation over $\ell$, followed by
\[
\sum_{r=0}^L\binom Lr^{-1}\le3,
\]
proves \eqref{eq:p27-time}.
\end{proof}

\begin{lemma}\label{s1}[Spatial convolution estimate]
For any $ k,n,n' \in\mathbb Z^d$, we have
\begin{equation}\label{eq:p27-space}
\begin{split}
&\sum_{k\in\mathbb Z^d}
e^{-\rho\phi(|n-k|)-\rho\phi(|k-n'|)}
h(n-k)h(k-n')\le
C_d e^{-\rho\phi(|n-n'|)}h(n-n').
\end{split}
\end{equation}
\end{lemma}
\begin{proof}
Monotonicity and subadditivity of $\phi$ give
\[
\phi(|n-n'|)
\le \phi(|n-k|)+\phi(|k-n'|).
\]
It therefore suffices to establish
\[
h*h\le C_dh.
\]
For $z,y\in\mathbb Z^d$, at least one of
\[
|y|\ge |z|/2,\qquad |z-y|\ge |z|/2
\]
holds. On the first region,
\[
h(y)\le 2^{d+2}h(z),
\]
and on the second region the same estimate holds for $h(z-y)$.
Since $h\in\ell^1(\mathbb Z^d)$,
\[
(h*h)(z)
\le 2^{d+3}\|h\|_{\ell^1}h(z).
\]
This proves \eqref{eq:p27-space}.
\end{proof}

\begin{lemma}[Radial weights after optimizing the smoothing order]\label{C}
For $T\in\mathbb{N}$ and $x\ge1$, define
\[
 F_T(x)=\sup_{m\in\mathbb{N}}
 \left\{\tau m\log x-\log\frac{W_{T+m}}{W_T}\right\},\qquad
 G_T(x)=W_Te^{-F_T(x)}=\inf_mW_{T+m}x^{-\tau m}.
\]
If $y\ge x$, $y/x\le10/9$, and $y-x\le d_*$, then, for any fixed $\rho>0$,
\begin{equation}\label{eq:transfer}
 G_T(x)\le C_{\sigma,\kappa,\tau,\rho}
 e^{\rho\phi(d_*)/2}G_T(y),
\end{equation}
with a constant independent of $T,x,y,d_*$.
\end{lemma}

\begin{proof}
The supremum is attained at a finite integer $m_y$. If $m_y\ge1$, comparison with $m_y-1$ gives
\[
 \kappa m_y^{p-1}\le\tau\log y,\qquad
 F_T(y)-F_T(x)\le\tau m_y\log(y/x)
 \le C(\log y)^{\sigma-1}\log(y/x).
\]
If $y\le(e+d_*)^2$, the last expression is at most $C[\log(e+d_*)]^{\sigma-1}$; otherwise,
$\log(y/x)\le (y-x)/x\le C d_*/y\le C y^{-1/2}$, so the expression is uniformly bounded.
Since $\sigma-1<\sigma$, both cases are controlled by $C+\rho\phi(d_*)/2$. The assertion is immediate if $m_y=0$.
\end{proof}

\begin{lemma}[Small-gap intercluster terms are smoothing]\label{D}
Let $0<\tau\le\tau_0$. For $r\ge1$ and
$S\in\widetilde{\mathcal{L}}_{\tau,\mathrm{ND}}^{-(r-1)}(\delta)$, $\delta\in\{0,1\}$,
let $S^{\mathrm{gap}}$ denote the matrix part satisfying the large-gap condition. Then
\[
 S-S^{\mathrm{gap}}\in\mathcal{R}_{r-1}^{-\infty},\qquad
 |S-S^{\mathrm{gap}}|_{r-1}\le C\|S\|_{r-1,\delta}.
\]
The parameters $M,\kappa$ must also be chosen initially large enough to absorb the following fixed moment estimate.
\end{lemma}

\begin{proof}
Intercluster separation and the small-gap condition give $|n-n'|\ge cX_{nn'}^{\tau_0}$. Adjusting constants at low frequencies yields
\[
 e^{-\rho\phi(|n-n'|)}\le C e^{-c_\rho(\log X_{nn'})^\sigma}
 \le C e^{b_gN^p}X_{nn'}^{-\tau N},
\]
where $c_\rho,b_g>0$ are fixed. The temporal index is at most $q+r-1$. Choose $\kappa\ge b_g$ and use
$W_{q+r-1+N}/W_{q+r-1}\ge e^{\kappa N^p}$. The original radial factor of negative order can simply be discarded.
\end{proof}

\end{document}